\documentclass[11pt]{amsart}

\newcommand{\ifn}[1]{} 

\usepackage{color, epsfig}

\usepackage{amssymb, amsmath, amscd, amsthm, graphicx, psfrag}
\usepackage[matrix,arrow,curve]{xy}
\xyoption{dvips}              

\newcommand{\calG}{{\mathcal{G}}}
\newcommand{\calO}{{\mathcal{O}}}
\newcommand{\calPo}{{\mathcal{P}_0}}

\newcommand{\calL}{{\mathcal{L}}}

\numberwithin{equation}{section}

\DeclareMathOperator*{\holim}{holim}
\DeclareMathOperator{\hoTot}{hoTot}
\DeclareMathOperator{\Tot}{Tot}
\DeclareMathOperator{\Hom}{Hom}

\newcommand{\BK}{{\mathbb K}}
\newcommand{\BR}{{\mathbb R}}

\newcommand{\BQ}{{\mathbb Q}}
\newcommand{\C}{{\operatorname{C}}}
\newcommand{\Ass}{{\mathcal{ASS}}}
\newcommand{\Poiss}{{\mathcal{POISS}}}
\newcommand{\COASS}{{\mathcal{ASS}}}
\newcommand{\unit}{{\mathbf{1}}}
\newcommand{\Emb}{{\operatorname{Emb}}}
\newcommand{\Imm}{{\operatorname{Imm}}}
\newcommand{\bEmb}{\overline{\operatorname{Emb}}}

\newcommand{\mult}{\mathrm{mult}}

\newcommand{\Kd}{\mathcal{K}_{d}}
\newcommand{\Kdn}{\mathcal{K}_{d[n]}}
\newcommand{\Ko}{{{\mathcal K}_1^{(0)}}}

\newcommand{\Nerve}{{\operatorname{N}}}
\newcommand{\norm}{{\operatorname{N}}}
\newcommand{\Ho}{{\operatorname{H}}}
\newcommand{\HH}{{\operatorname{HH}}}
\newcommand{\Dd}{{{\mathcal{D}}_d}}
\newcommand{\Top}{\mathrm{Top}}
\newcommand{\Ab}{\mathrm{Ab}}
\newcommand{\Ch}{\mathrm{Ch}}
\newcommand{\CDGA}{\mathrm{CDGA}}
\newcommand{\DGL}{\mathrm{DGL}}
\newcommand{\Apl}{A_{PL}}
\newcommand{\OmPA}{\Omega_{PA}}
\newcommand{\Fd}{{{\mathcal F}_d}}
\newcommand{\Fo}{{{\mathcal F}_1^{(0)}}}
\newcommand{\fanicop}[2]{\widehat{{#1}}_{\{{#2}\}}}

\newcommand{\fanic}{fanic}

\newcommand{\bead}{bead}
\newcommand{\calC}{{\mathcal{C}}}
\newcommand{\calF}{{\mathcal{F}}}
\newcommand{\id}{\mathrm{id}}

\newcommand{\calR}{{\mathcal{R}}}
\newcommand{\calM}{{\mathcal{M}}}
\newcommand{\calK}{{\mathcal{K}}}
\newcommand{\calKd}{{\mathcal{K}_d}}

\newcommand{\nso}{nonsymmetric operad}
\newcommand{\quism}{\stackrel{\simeq}{\to}}

\theoremstyle{plain}
\newtheorem{thm}{Theorem}[section]

\newtheorem{prop}[thm]{Proposition}

\newtheorem{lemma}[thm]{Lemma}
\newtheorem{cor}[thm]{Corollary}

\theoremstyle{definition}
\newtheorem{definition}[thm]{Definition}
\newtheorem{eg}[thm]{Example}
\newtheorem{rem}[thm]{Remark}

\newcommand{\refT}[1]{Theorem~\ref{T:#1}}
\newcommand{\refC}[1]{Corollary~\ref{C:#1}}
\newcommand{\refP}[1]{Proposition~\ref{P:#1}}

\newcommand{\refL}[1]{Lemma~\ref{L:#1}}
\newcommand{\refE}[1]{Equation~$(\ref{E:#1})$}
\newcommand{\refS}[1]{Section~\ref{S:#1}}
\newcommand{\refZ}[1]{$(\ref{#1})$}

\begin{document}


\title{The rational homology of spaces of long knots in codimension $>$2}


\author{Pascal Lambrechts}
\address{Universit\'e Catholique de Louvain, 2 Chemin du Cyclotron, B-1348 Louvain-la-Neuve, Belgium}
\email{lambrechts@math.ucl.ac.be}
\urladdr{http://milnor.math.ucl.ac.be/plwiki}
\author{Victor Turchin}
\address{Universit\' Catholique de Louvain, 2 Chemin du Cyclotron, B-1348 Louvain-la-Neuve, Belgium.\\
University of Oregon, USA.\\ 
Institut des Hautes Etudes Scientifiques, France.}
\email{turchin@math.ucl.ac.be} \urladdr{http://www.math.ucl.ac.be/membres/turchin/}
\author{Ismar Voli\'c}
\address{Department of Mathematics, Wellesley College,
Wellesley, MA}
\email{ivolic@wellesley.edu}
\urladdr{http://palmer.wellesley.edu/\~{}ivolic}

\thanks{The first author is Chercheur Qualifi\'e at F.N.R.S.
The second author was supported in part by grants NSH-1972.2003.01 and RFBR 05-01-01012a. 
 The third author was supported in part by the National Science Foundation grant DMS 0504390.  }


\begin{abstract}
We determine the rational homology of the space of long knots in $\BR^d$
for $d\geq4$.
Our main result is that the Vassiliev spectral sequence computing this rational homology collapses
at the $E^1$ page. 
As a corollary we get that the homology of long knots (modulo immersions) is the Hochschild homology of the Poisson algebras operad with a bracket
 of degree $d-1$, which can
be obtained as the homology of an explicit graph complex and is in theory completely computable.

Our proof is a combination of a relative version of Kontsevich's formality of the little $d$-disks  operad and of
Sinha's cosimplicial model for the space of long knots arising from Goodwillie-Weiss embedding calculus.
 As another ingredient in our proof, we
 introduce a generalization of a construction 
 that associates a cosimplicial
object to a multiplicative operad. Along the way we also establish some results 
about the Bousfield-Kan spectral sequences of a truncated cosimplicial space.
\end{abstract}
 
\maketitle


\sloppy

\ifn{For some reason, I can't put footnotes right into the abstract, but here is what I changed:
\\1) I replaced ``$E^*$-page" by ``$E^*$ page" everywhere in the paper.  The first looks awkward.  Not only that, but there was an ``$n$th-term" somewhere later in the paper which is definitely not right.  It should be ``$n$th term".
\\2) There was an equation in the abstract which I took out.  It's bad form to put that much math notation in an abstract.
\\3) We really need to get into the habit of saying ``Goodwilie-Weiss embedding calculus" and not just ``Goodwillie embedding calculus".
\\4) I took out ``construction of McClure and Smith which..."  The reason is that Jim actually told me once not to say this.  This construction is due to Gerstenhaber and Voronov, and he wants the credit to go to them.  I've tried to change other things in the paper to accomodate this a little bit.
\\5) ``In the course we also establish" is not correct.
}


\section{Introduction}



A  {\em long knot} is a smooth embedding $f\colon\BR\hookrightarrow\BR^d$, 
$d\geq 3$, that coincides with a fixed  linear embedding 
$\epsilon\colon\BR\hookrightarrow\BR^d$ 
outside a compact set.
We denote by $\Emb(\BR,\BR^d)$
the space of all long knots equipped with the weak ${\mathcal{C}}^\infty$ topology.
This space, and in particular its homology,  has been under investigation for more than twenty years.
One of the first important tools in this study was a spectral sequence for computing the homology of spaces of knots constructed  by\ifn{``Invented by..." is awkward.  You invent something physical, like a machine, but not really a concept.} Vassiliev in the late eighties \cite{Vass:CKS}.  This spectral sequence sparked a lot of interest,
especially in dimension $d=3$, where it is related to the theory of finite type knot invariants.

Independently, Goodwillie and Weiss suggested another approach for the study of knots and more general
spaces of embeddings which is now known as \emph{embedding calculus}. In particular, it is suggested in \cite[Example 5.1.4]{GKW:survey} that this approach should also give a spectral sequence 
for computing the homology of spaces of long knots.  Indeed, it later turned out that it does, and that this spectral sequence was
equivalent to Vassiliev's.
Goodwillie-Weiss embedding calculus for knots was developed \ifn{``Developped" is a common mistake; it should be ``developed".} further by Sinha in \cite{S:TSK},
who also emphasized the connection with the little $d$-disks operad in \cite{S:OKS}.

Vassiliev spectral sequence arises from a study of the \emph{discriminant set}, i.e. the complement of the set of knots 
in the space of all smooth maps $\BR\to\BR^d$ with fixed behavior\ifn{``Behaviour" is British.} at infinity
(or more precisely in some finite-dimensional approximation of that space). 
In other words, the discriminant set is the set of maps with singularities.  This set admits a nice stratification, which in turn yields a natural filtration from which the spectral sequence is constructed. 
A classification of the singularities gives the $E^1$ page of this spectral sequence a combinatorial description in terms
of certain graphs such as chord diagrams familiar from finite type knot theory.
 When $d\geq4$, this spectral sequence converges to
the homology of the space of long knots~\cite[Section~6.6]{Vass:CKS}.

Vassiliev has conjectured a stable splitting of the resolved discriminant which would imply that his spectral sequence collapses at the $E^1$ page.  This collapse was proved rationally by Kontsevich in dimension $d=3$ 
along the diagonal ${E^1_{-p,p}}$. 
The proof uses the famous \emph{Kontsevich Integral}, a map that realizes\ifn{``Realises" is British.} all finite type invariants \cite{K:Vass}.
Kontsevich further claimed in \cite[Theorem 2.3]{K:Fey} that his integration approach can
be generalized  for $d\geq4$ to give a proof of 
the  collapse\ifn{``collapsing of the SS" sounds awkward.  I've changed it to ``collapse" everywhere.} of the rational Vassiliev spectral sequence everywhere.
 In  \cite{CCRL},
Cattaneo, Cotta-Ramusino, and Longoni filled in some details of that program and proved
the collapse along the main diagonal.
As far as we know, however, no complete proof of the rational collapse has yet appeared.

In this paper, we give a proof of Vassiliev's conjecture over the rationals.
\begin{thm}\label{T:VassConj}
For  $d\geq4$,
the Vassiliev spectral sequence computing the rational homology of the space of long knots
$\Emb(\BR,\BR^d)$ collapses at the $E^1$ page.
\end{thm}

To prove this theorem, we take a very different approach than that of Kontsevich, namely the Goodwillie-Weiss calculus of the embedding functor and Sinha's cosimplicial model for spaces of knots arising from this theory. 
Before explaining this further, it is convenient to introduce
 a variation\ifn{There is a lot of ``some" as in ``there is some SS" or ``to introduce some variation" etc.  This sounds awkward.  It's better to use `a", like ``to introduce a variation".} of the space of long knots. Consider first the  \emph{space of long immersions}
 $$\Imm(\BR,\BR^d):=\{f\colon\BR\looparrowright \BR^d:f \textrm{ is an immersion that coincides with } \epsilon
 \textrm{ outside a compact set}\}.$$
There is an inclusion $\iota\colon\Emb(\BR,\BR^d)\hookrightarrow\Imm(\BR,\BR^d)$ and its 
homotopy fiber\ifn{``Fibre" is British.}, $\bEmb(\BR,\BR^d)$, is called the \emph{space of long knots modulo immersions}.
By Smale-Hirsh theory \cite{Smale:ClassImmSn}, there is a weak equivalence 
$\Imm(\BR,\BR^d)\simeq\Omega S^{d-1}$.
 Moreover $\iota$ is
 null-homotopic  \cite[Proposition 5.17]{S:OKS}, so we have a weak equivalence
 \begin{equation}\label{E:bEmbsplit}
 \bEmb(\BR,\BR^d)\simeq\Emb(\BR,\BR^d)\times\Omega^2S^{d-1}.
 \end{equation}

 In \cite{S:OKS} Sinha constructs a cosimplicial space
$$
\xymatrix{
\Kd^\bullet = 
\big(\Kd(0) \ar@<0.7ex>[r] \ar@<-0.7ex>[r] &
\Kd(1) \ar[l] \ar[r] \ar@<1.2ex>[r]  \ar@<-1.2ex>[r] &
\Kd(2) \ar@<0.6ex>[l]  \ar@<-0.6ex>[l]
\cdots \big),
}
$$
where $\Kd=\{\Kd(n)\}_{n\geq0}$ is a topological operad homotopy
equivalent to the little $d$-disks operad, called
 the \emph{Kontsevich operad}. This operad turns out to be \emph{multiplicative}, i.e. there exists a map $\{*\}_{n\geq0}\to\Kd$ from the nonsymmetric associative topological operad which consists of 
the one-point space in each degree, to the Kontsevich operad. 
The cofaces and codegeneracies of $\Kd^\bullet$
 are induced from this multiplicative structure via a general construction 
 of Gerstenhaber and Voronov \cite{GV:HGMSO} which we
recall in \refS{cos-multop} using the McClure-Smith point of view \cite[Section 3]{MS:Del}. The main result of \cite{S:OKS} is that for $d\geq 4$ 
the homotopy totalization\ifn{``totalisation" is British.} of that cosimplicial space, $\hoTot(\Kd^\bullet)$,  is weakly equivalent to $\bEmb(\BR,\BR^d)$.

The homology Bousfield-Kan spectral sequence \cite{B:HSS} associated to $\Kd^\bullet$ converges
to $\Ho_*(\hoTot(\Kd^\bullet))$ when $d\geq4$ \cite[Theorem 7.2]{S:OKS}. 
Our main result is 
\begin{thm}\label{T:BKSinha-coll}
For $d\geq4$ the homology Bousfield-Kan spectral sequence associated to Sinha's cosimplicial space $\Kd^\bullet$ collapses
at the $E^2$ page rationally.
\end{thm}
Before looking at\ifn{``Looking to" is incorrect.  ``Looking at" is correct.} the consequences of this theorem, we give an overview of its proof.
The key idea is a relative version of Kontsevich's theorem on the formality of
 the little $d$-disks operad \cite[Section 3]{K:OMDQ}, which states that there is a chain of quasi-isomorphisms of operads between the singular chains 
of the little $d$-disks operad and its homology with real coefficients.
 This is therefore true for any operad weakly equivalent to the little $d$-disks operad.  In particular there is a chain of quasi-isomorphisms $\C_*(\Kd;\BR)\simeq\Ho_*(\Kd;\BR)$.  If one could deduce from this the formality of the cosimplicial space
$\Kd^\bullet$, the collapse of the homology Bousfield-Kan  spectral sequence would immediately follow from \refP{formalcollapse}. For this it would be enough to know, by McClure-Smith, that $\Kd$ is formal as a \emph{multiplicative} operad. But this does not appear to be easy to establish (see Remark \ref{R:nonmult}), and 
we do not know whether $\Kd^\bullet$  is formal.\ifn{There used to be a ``whether the cosimplicial space $\Kd^\bullet$ is formal" here.  But this is redundant, because we've been saying all along that $\Kd^\bullet$ is a cosimplicial space, and the reader at this point doesn't need to be reminded of it.  In general, there is a lot of this is the paper -- we name the thing and then give its notation, as in ``Here is a map $f$.  We use the map $f$ to do this.  Notice that the map $f$ does this".  In the last two sentences, ``the map" can and should be removed.}
 
This is why we make a detour\ifn{As mentioned in an earlier footnote, there was a ``some" issue here before.  It used to say ``we will make some detour...".  which is very awkward.} through the \emph{Fulton-MacPherson operad} $\Fd$. This operad  is 
homotopy equivalent to $\Kd$ but is more suitable for proving formality results as those in  \cite{K:OMDQ}.
It is not multiplicative, but it is ``multiplicative up to homotopy'' in the sense that it admits a map from the contractible $A_\infty$ Stasheff associahedral
 operad. The key formality statement that we prove is that when $d\geq 3$, $\Fd$ is
formal as an operad which is  ``multiplicative up to homotopy'' (\refL{collinearities}). However, the problem is that we cannot directly 
 associate to $\Fd$ a cosimplicial space in the spirit of McClure-Smith.\ifn{This used to say ``\`a la McClure-Smith" which is used in English and there's no real problem with it, but it nevertheless sounds very pretentious.}  
We instead construct certain finite diagrams $\fanicop{\Fd}{n}$ of spaces 
which we call \emph{fanic diagrams}. These are built out of the up-to-homotopy multiplicative structure 
on $\Fd$ and are in a sense 
a rigid version of a ``cosimplicial space up to homotopy'' analogous to $\calK^\bullet_d$, or more precisely analogous to its $n$th truncation.
Our formality statement implies that $\fanicop{\Fd}{n}$ is formal (\refT{FormalityFanic}) and hence the homology Bousfield-Kan  spectral sequence of the 
cosimplicial replacement of that diagram collapses at the $E^2$ page.  Thus the same must be true for the cosimplicial replacement of
 the $n$th truncation $\Kdn$ of $\Kd^\bullet$. Using the fact that, for each $n\geq0$, the homology spectral sequence of $\Kdn$ 
collapses, along with a strong convergence result for
 Bousfield-Kan spectral sequences as it applies to $\Kd^\bullet$,
we deduce \refT{BKSinha-coll}.\ifn{There was a lot of tense mixing in the past few sentences, of the sort ``We will do this, which implies that."  The first part of this sentence is in the future tense, and the second in the present, which makes the sentence incorrect.}

\vskip 6pt
\noindent
We now come back\ifn{This used to say ``We come back to the..." but that sounds awkward.  And it really kind of means that we regularly come back to this, or we come back to it every once in a while.  I.e. ``we come back" means that there is a continuous action.} to the Vassiliev spectral sequence and explain its link with the Bousfield-Kan spectral sequence.  In his thesis, the second author found a more conceptual description of the $E^1$ page of the Vassiliev 
spectral sequence than the original 
combinatorial one arising from the classification of
singularities in the discriminant set. To explain this,  
recall that the \emph{Poisson operad with bracket of degree $d-1$}, $\Poiss_{d-1}$,
 is the operad  encoding Poisson algebras, i.e. graded commutative algebras equipped with a Lie bracket
 of degree $(d-1)$ which is a graded  derivation with respect to the multiplication (see, for instance, Example (d) in \cite[Section 1]{T:HSLK}).
The work of \ifn{This used to say ``Actually, the work of...", but ``actually" is unnecessary and awkward.  It implies that something slightly unexpected and connected in an unusual way to what it follows is going to be said.  But that's not the case here.  I took out a lot of ``actually"s throughout the paper.} Fred Cohen \cite{CLM:HILS} implies that, for $d\geq2$, $\Poiss_{d-1}$ is just the homology of the topological little $d$-disks operad.
 Since  Poisson algebras admit an associative  multiplication $m\in\Poiss_{d-1}(2)$, we can define a differential
 $$\delta\colon\Poiss_{d-1}(n)\longrightarrow\Poiss_{d-1}(n+1)$$
 where $\delta=[m,-]$ is  the Gerstenhaber bracket and $\delta^2=0$ by associativity\ifn{There used to be a comma before the reference here.  It's not needed.} \cite[Section 3]{T:HSLK}.
 The homology of the cochain complex $(\Poiss_{d-1}(*),\delta)$ is called the \emph{Hochschild homology} of that
 multiplicative operad, denoted by \ifn{Again an issue with comma.  I added ``denoted by".  Without it, the sentence is awkward.}$\HH_*(\Poiss_{d-1})$. 
 This complex is in fact the deformation complex of the morphism of operads $\Ass\to\Poiss_{d-1}$, where $\Ass$ is the associative algebras operad~\cite{KS:DAODC}.
  Finally, a slight variation of Vassiliev's spectral sequence for $\Emb(\BR,\BR^d)$ produces a spectral sequence for $\bEmb(\BR,\BR^d)$ and
we have
 \begin{thm}[Turchin]\label{T:TurchinE1}
 The $E^1$ page of the Vassiliev spectral sequence computing $\Ho_*(\bEmb(\BR,\BR^d))$
 is isomorphic to the Hochschild homology of the Poisson operad with bracket of degree $d-1$.
\end{thm}
This statement was first established for the Vassiliev spectral sequence for $\Emb(\BR,\BR^d)$ \cite{T:HSLK},
and the above analogous result is a combination of \cite[Theorem 8.4]{T:OSCD-x2} and 
\cite[Proposition 3.1 and Lemma 4.3]{T:W1THH}.

On the other hand, the $E^2$ page of the Bousfield-Kan spectral sequence for $\Kd^\bullet$ is by definition the homology
of the conormalisation of the cosimplicial abelian group $\Ho_*(\Kd^\bullet)$. Since the cofaces are induced by the multiplicative
structure on $\Kd$ and since the homology of that multiplicative topological operad is exactly the
multiplicative Poisson operad, it is easy to deduce that
this $E^2$ page is also isomorphic to  $\HH_*(\Poiss_{d-1})$ \cite[Corollary 1.3]{S:OKS}.
This implies that
 the $E^1$ page of the Vassiliev spectral sequence for $\bEmb(\BR,\BR^d)$  and the 
$E^2$ page of the homology spectral sequence for $\Kd^\bullet$  are isomorphic.  It can be shown that this isomorphism is just a regrading \cite[Proposition 0.1]{T:OSCD-x2}. 
Since both of these spectral sequences converge to $\Ho_*(\bEmb(\BR,\BR^d);\BQ)$  for $d\geq4$,
\refT{BKSinha-coll} implies that the Vassiliev spectral sequence collapses at $E^1$.
As an immediate consequence of this and
 \refZ{E:bEmbsplit}  we have the following
\begin{cor}\label{C:HEmb} For $d\geq 4$,
$$
\Ho_*(\Emb(\BR,\BR^d);\BQ)\otimes \Ho_*(\Omega^2 S^{d-1};\BQ)\cong \Ho_*(\bEmb(\BR,\BR^d);\BQ)\cong\HH_*(\Poiss_{d-1}).
$$
\end{cor}
The rational homology of $\Omega^2 S^{d-1}$ is isomorphic to a free graded commutative algebra on one or two generators depending on the parity of $d$, so it is very simple.
In addition\ifn{This used to say ``on the other hand", but that's misleading because the reader is expecting something different, almost contradictory to what was just said, but what follows agrees with what was said.  Namely, we say homology of spheres is easy.  In addition, homology of configuration spaces
is easy.  Not ``On the other hand, homology of configurations is hard".}, the $E^2$ page of the Bousfield-Kan spectral sequence for $\Kd^\bullet$,
or equivalently the Hochschild homology of the Poisson operad,
can be expressed as the homology of an explicit cochain complex described in terms of the homology of configuration spaces in $\BR^d$ which are well understood.  
In particular, there is an algorithm for computing this homology which has been
used in~\cite{T:OSCD-x2}.  
Thus, in theory, $\Ho_*(\Emb(\BR,\BR^d);\BQ)$ is completely computable. However, it appears that the 
algorithmic complexity is exponential and the only feasible computations are all still in low degrees.  What is necessary, and is still lacking, is a deeper understanding of the structure of $\HH_*(\Poiss_{d-1})$.

It is easy to see that the homology of the cosimplicial space $\Kd^\bullet$ depends, up to an obvious regrading, only on the parity of $d$.
This implies that, up to regrading, all the $E^2$ pages of the homology Bousfield-Kan  spectral sequences 
of $\Kd^\bullet$ 
for $d$ even (respectively $d$ odd) are isomorphic.
 This completes the proof of the characterization of the rational homology of spaces of long knots as stated in \cite[Theorem 2.3]{K:Fey}.
Furthermore, the operation of stacking of long knots gives a multiplication which makes $\Emb(\BR,\BR^d)$ into an $H$-space.  Therefore the rational homotopy type of
$\Emb(\BR,\BR^d)$ for $d\geq4$ is completely determined by its rational homology and is thus
 virtually known by our main theorem.

Coming back to the original Vassiliev conjecture, notice that by \cite[Theorem 5.1]{T:W1THH} and results of
\cite[Section~IV]{Vass:inv-discr}, $E^1$ page of the Vassiliev spectral
sequence for $\Emb(\BR,\BR^d)$ and the corresponding $\overline{E^1}$ page of the Vassiliev spectral
sequence for $\bEmb(\BR,\BR^d)$ are related by an isomorphism
$E^1\otimes \Ho_*(\Omega^2 S^{d-1})\cong\overline{E^1}$.
Therefore \refC{HEmb} and the isomorphism $\overline{E^1}\cong \HH_*(\Poiss_{d-1})$ from \refT{TurchinE1}
imply the collapse at $E^1$ of the classical Vassiliev spectral sequence for  $\Emb(\BR,\BR^d)$, which is 
our promised \refT{VassConj}.

\vskip 6pt
\noindent
The \emph{homotopy} Bousfield-Kan spectral sequence computing the homotopy groups of $\Emb(\BR,\BR^d)$
from the cosimplicial abelian group $\pi_*(\Kd^\bullet)$ has been studied elsewhere \cite{ALTV:coformal-x2}. 
Using the remarkable fact that the little $d$-disks operad is not only formal
(i.e. rationally determined by its homology), but also \emph{coformal} (i.e. rationally determined by its homotopy), it has been shown 
in \cite{ALTV:coformal-x2} that \refT{BKSinha-coll} implies that the homotopy spectral sequence tensored with the rationals  also collapses at the $E^2$ page.  A strong connection between the $E^2$ pages of the rational homology and the homotopy spectral sequences is exhibited in \cite{LT:piGC-x1} where it is also shown that the homotopy groups of the space of long knots are also the homology of a certain graph complex, smaller
than the one used for computing the homology of $\Emb(\BR,\BR^d)$. Some computations in this $E^2$ page appeared in \cite{ScaS}.

A natural next question is whether our approach gives the collapse of the Vassiliev spectral sequence for long knots in $\BR^3$.  There are many difficult issues here, not the least of which is that it is not even clear what the Vassiliev spectral sequence 
converges to.  However, we are nevertheless able to prove a certain collapsing result. Indeed, \refP{Kdncoll} states 
that the homology Bousfield-Kan  spectral sequence
of the cosimplicial replacement of the $n$th truncation $\Kdn$ collapses at $E^2$. This is true for $d\geq 3$,
but to deduce the collapse of the Bousfield-Kan spectral sequence for $\Kd^\bullet$ itself we need some convergence hypotheses
which are only true for $d\geq 4$. It at least seems plausible that this spectral sequence collapses for $d=3$.  Another interesting question is whether our result can be extended over the integers.

Finally, it is possible to extend our results to spaces of embeddings of more general manifolds than $\BR$ in $\BR^d$.  This is done by generalizing the following slogan which summarizes our proof of \refT{BKSinha-coll}:

\vskip 8pt

\begin{tabular}{ccc}
Goodwillie-Weiss embedding calculus&& collapse of\\
 {\large+} &$\Longrightarrow$& spectral sequences\\
Kontsevich formality of the little $d$-disks  operad &&for spaces of embeddings
\end{tabular}

\vskip 8pt

\noindent This slogan was taken much further in \cite{ALV:HQEmb-x1}.
The authors of that paper use Weiss' orthogonal calculus to prove a collapse result for a spectral sequence computing the rational homology of $\bEmb(M,\BR^d)$
 for any compact manifold $M$ and $d$ large enough. In particular, it is shown that for large $d$, the rational homology of the space of embeddings modulo immersions,
 $\Ho_*(\bEmb(M,\BR^d);\BQ)$, depends only on the rational homotopy type of $M$. 

\subsection{Organization of the paper}
In \refS{KdFdKdbullet} we review basic facts on cosimplicial spaces and diagrams.  We also recall Sinha's cosimplicial model for the space of long knots, the notion of a multiplicative
operad and the associated Gerstenhaber-Voronov/McClure-Smith cosimplicial object, and the definitions of the Kontsevich and Fulton-MacPherson operads.
Formality, which comes from rational homotopy theory, is recalled in \refS{formalityBKSS} where we also prove that the homology
Bousfield-Kan spectral sequence of a formal cosimplicial space, or more generally of the cosimplicial replacement of a formal
finite diagram, collapses at the $E^2$ page. In \refS{SStotn} we study a condition on a cosimplicial space that guarantees that the
Bousfield-Kan spectral sequences for its truncations converge. In \refS{Fans} we introduce the category of \emph{fans}, which
is a variation of the truncated cosimplicial category $\Delta[n]$, and in Section \ref{S:Fanicnso} show how to associate to any morphism of
nonsymmetric operads\ifn{This used to say ``morphism of nonsymmetric operad", which is incorrect -- a morphism involves 2 operads, so we need a plural.} a diagram shaped on the category of fans.  This generalizes the Gerstenhaber-Voronov/McClure-Smith cosimplicial diagram associated to 
a multiplicative operad. In \refS{formalfanicFd} we prove a relative version of Kontsevich's formality of the little $d$-disks operad 
and deduce the formality of the fanic diagrams associated to the Kontsevich multiplicative operad. In the last section we collect these
results to give a proof of our main theorem. 
\subsection{Acknowledgments}

We are very grateful to Greg Arone for many conversations and encouragement, as well as for arranging
several visits of the first author to University of Virginia.
We also thank  Peter Bousfield, Emmanuel Dror-Farjoun, and
 Bill Dwyer  for answering   our many questions on
 the Bousfield-Kan spectral sequences; Tom Goodwillie and Dev Sinha for fruitful conversations; and 
  Ryan Budney for help with figures.


\section{The Kontsevich and Fulton-MacPherson operads and
 Sinha's cosimplicial model for the space of long knots}\label{S:KdFdKdbullet}

In this section we review Sinha's cosimplicial model for the space of long knots and its relation to the Kontsevich and Fulton-MacPherson operads. To start, we review some standard facts about cosimplicial spaces and diagrams as well as the construction of a cosimplicial space associated to a multiplicative operad. 

 \subsection{Cosimplicial spaces, homotopy totalizations, Bousfield-Kan spectral sequences, cosimplicial replacement of diagrams, and
 left cofinal functors.}\label{S:cos-hotot}
The standard references for the following basic terminology and facts about cosimplicial objects are \cite[X]{BK:LNM} and~\cite[Chapter~8]{Weibel}. 

The simplicial category  $\Delta$ has ordered sets $[n]:=\{0,1,2,..., n\}$, $n\geq 0$, as objects and order-preserving maps as morphisms.
 All morphisms in $\Delta$ are compositions of \emph{cofaces} $d^i\colon[n]\to[n+1]$  and \emph{codegeneracies}
 $s^j\colon[n]\to[n-1]$. A \emph{cosimplicial object} in a category $\calC$ is a covariant
 functor from $\Delta$ to $\calC$. Dually, a \emph{simplicial object} in $\calC$ is a contravariant
 functor from $\Delta$ to $\calC$.
In particular  a \emph{(pointed) cosimplicial space} is a covariant functor 
$$X^\bullet\colon\Delta\longrightarrow\Top\quad(\textrm{or }\Top_*) $$ from the simplicial category
to the category of (pointed) spaces.
The \emph{standard cosimplicial space} is the cosimplicial space $\Delta^\bullet$ where $\Delta^n$ is the standard geometric $n$-simplex, 
the cofaces are defined from the inclusions of faces $\Delta^n\hookrightarrow\Delta^{n+1}$, and the codegeneracies 
are suitable affine projections $\Delta^n\to\Delta^{n-1}$.

The \emph{totalization} of a cosimplicial space $X^\bullet$, denoted by $\Tot X^\bullet$, 
is the space of natural maps from the standard  cosimplicial space $\Delta^{\bullet}$  to $X^\bullet$,
 $$
 \Tot X^\bullet:=\Hom_\Delta(\Delta^\bullet,X^\bullet).
 $$  
When the cosimplicial space $X^\bullet$ is \emph{fibrant} \cite[X.4.6]{BK:LNM}, $\Tot(X^\bullet)$ is homotopy equivalent to the
homotopy limit $\holim_{\Delta}X^\bullet$ \cite[XI.4.4]{BK:LNM}. As explained in \cite[Section 2.7]{B:HSS},  the techniques
 of \cite[pp.279-280]{BK:LNM} imply that any cosimplicial space $X^\bullet$  admits
a weakly equivalent fibrant functorial replacement $\tilde X^\bullet$.  One then defines 
  the \emph{homotopy totalization}
of $X^\bullet$ by
$$\hoTot X^\bullet:=\Tot \tilde X^\bullet.$$
Since by definition $\tilde X^\bullet$ is fibrant and  weakly equivalent to $X^\bullet$, $\hoTot X^\bullet$ is always weakly equivalent
to $\holim_{\Delta}X^\bullet$.  This homotopy totalization $\hoTot$ is also weakly
 equivalent to the homotopy invariant totalization $\widetilde{\Tot}$ used 
in \cite{S:OKS}.

Consider the full subcategory $\Delta[n]\subset\Delta$ consisting of objects $[0],\ldots,[n]$. The \emph{$n$th truncation}
of a cosimplicial object $X^\bullet$ in $\calC$ is the composite
 $$X_{[n]}:\Delta[n]\hookrightarrow\Delta\stackrel{X^\bullet}\longrightarrow\calC.$$
The \emph{$n$th partial homotopy totalization} of $X^\bullet$ is defined as
$$\hoTot^n(X^\bullet):=\holim_{\Delta[n]}X_{[n]}.$$

To any cosimplicial space $X^\bullet$ one can associate a second quadrant homology  Bousfield-Kan spectral sequence
with cofficients in an abelian group $A$ \cite{B:HSS}, and a homotopy Bousfield-Kan spectral sequence if $X^\bullet$ is pointed \cite{BK:LNM}.  These converge under favorable circumstances to $\Ho_*(\hoTot X^\bullet;A)$ or $\pi_*(\hoTot X^\bullet)$.  We state certain strong convergence conditions in \refS{SStotn}.

A \emph{diagram of spaces} is a covariant functor  $F\colon I\to\Top$ where $I$ is a small category  called the \emph{shape} of the diagram.
We will also use the notion of a \emph{cosimplicial replacement of the diagram $F$} as defined in \cite[XI.5]{BK:LNM}, denoted by $\Pi^\bullet F$. 
The $n$th term of this cosimplicial space is given by 
$$
\Pi^n F:=\prod_{u=(i_0\stackrel{f_1}\leftarrow\cdots\stackrel{f_n}\leftarrow i_n)\in\Nerve_n(I)} 
F(i_0)
$$
where the simplicial set $\Nerve_\bullet(I)$ is the nerve of the category $I$. The cofaces and codegeneracies are obtained as suitable diagonal maps and projections.
By \cite[XI.5.2]{BK:LNM}, the
(homotopy) totalization of this cosimplicial space is weakly equivalent to the homotopy limit of $F$, i.e.
$$
\hoTot(\Pi^\bullet F)\simeq\holim_I F.
$$ 

If $\theta\colon J\to I$ is a functor then the $I$-diagram $F$ induces a $J$-diagram $\theta^*F:=F\circ \theta$. 
This \emph{change of shapes} functor $\theta$ is said to be \emph{left cofinal} if for every object $i$ in $I$,
the overcategory $\theta\downarrow i$ is contractible (see \cite[XI.9]{BK:LNM} with $\theta\downarrow i$ defined and denoted by $\theta/i$ in \cite[XI.2.2]{BK:LNM}). In this case there is a weak equivalence 
$$
\holim_I F\simeq\holim_J\theta^* F.
$$

\subsection{The construction of a cosimplicial object associated to a multiplicative operad following McClure and Smith}\label{S:cos-multop}
Here we recall the notion of the cosimplicial object associated to a multiplicative operad from \cite[Section 3]{MS:Del}.

 Let  $(\calC,\otimes,\unit)$ be a symmetric monoidal category where the object $\unit$ is the unit for $\otimes$. 
 A \emph{\nso}  $\calO=\{\calO(n)\}_{n\geq0}$ is a collection of objects of $\calC$ with all the properties of an operad except those having to do with the actions of the symmetric
groups \cite[Definition 3.12]{May:GILS}.  Notice that as part of the definition there is a unit morphism 
$\id\colon\unit\to\calO(1)$ and we also suppose that a (nonsymmetric) operad has an object in degree $0$, contrary to the
definition in \cite{MSS:operads}. We denote by 
$$
\circ_i\colon\calO(p)\otimes\calO(q)\longrightarrow\calO(p+q-1), \ \ \ 1\leq i\leq p,
$$ 
the usual insertion operations induced by the operad structure
  as defined for example in \cite[p. 7]{MSS:operads}.
The \emph{associative \nso}  $\Ass$ is  defined by $\Ass(n)=\unit$ for each $n\geq0$, with operadic structure maps the standard isomorphisms
$\unit^{\otimes(1+n)}\cong\unit$. 

 A \emph{multiplicative operad} is a \nso\ $\calO$ equipped with a morphism of \nso s
$\mu\colon\Ass\to\calO$, (\cite{GV:HGMSO},\cite[Definition 3.1 and Remark 3.2 (i)]{MS:Del}). Such a multiplicative structure on $\calO$
is equivalent to having morphisms $e\colon\unit\to\calO(0)$ and $m\colon\unit\to\calO(2)$ satisfying 
\begin{equation}\label{E:mult}
m\circ_1 m=m\circ_2m\quad \text{and}\quad m\circ_1 e= m\circ_2 e=\id.
\end{equation}
One can associate a cosimplicial object $\calO^\bullet$ to any multiplicative operad $\calO$ \cite[Section 3]{MS:Del}  
by defining $\calO^n=\calO(n)$ with the cofaces and codegeneracies
given by the following formulas.  For $x\in\calO(n)$,
\begin{align*}
d^0(x)&=m\circ_2x,\\
d^i(x)&=x\circ_i m,\quad\textrm{ for }1\leq i\leq n,\\
d^{n+1}(x)&=m\circ_1 x,\\
s^j(x)&=x\circ_j e,\quad\textrm{ for }1\leq j\leq n.
\end{align*}
It is easy to check that the cosimplicial identities are consequences of \refE{mult}.

There are also obvious dual notions of a \emph{cooperad}, of a \emph{coassociative cooperad} (also denoted by $\Ass$),
of a \emph{comultiplicative cooperad}, and of a simplicial object associated to a comultiplicative cooperad.


\subsection{The Kontsevich multiplicative operad and Sinha's cosimplicial model}\label{S:Kd}


Fix $d\geq 1$ and recall that a linear embedding  $\epsilon\colon\BR\hookrightarrow\BR^d$ has also been fixed.

The Kontsevich operad $\Kd=\{\Kd(n)\}_{n\geq0}$ is defined and studied in
\cite[Definition 4.1 and Theorem 4.5]{S:OKS}.
Each space $\Kd(n)$ is obtained as a suitable compactification of the ordered configuration space of $n$
points in $\BR^d$ modulo the action of $\BR\ltimes\BR^d$ by scaling and translation.  Another feature is that colinear configurations are identified. More precisely, $\Kd(n)$ is the closure of the image of the map
\begin{equation}\label{E:defKd}
\alpha_*=(\alpha_{ij})_{1\leq i<j\leq n}\colon C(n,\BR^d)\longrightarrow\prod_{1\leq i<j\leq n}S^{d-1}
\end{equation}
where $C(n,\BR^d)$ is the space of configurations of $n$ points in $\BR^d$ 
and $\alpha_{ij}\colon C(n,\BR^d)\to S^{d-1}$ is defined by
$\alpha_{ij}(x_1,\ldots,x_n)=(x_i-x_j)/\|x_i-x_j\|$.
It can be shown that this operad is homotopy equivalent to the classical little $d$-disks operad.

The spaces $\calK_1(n)$ turn out to be  homeomorphic to the discrete symmetric group  $\Sigma_n$ on $n$ letters, 
since all configuration on $\BR$ are colinear. 
Let $\Ko(n)$ be the path-connected component of $\calK_1(n)$ corresponding to the linearly ordered  configuration $(1,\ldots,n)$ 
on the line. Since  $\Ko(n)$ is a one-point space, it is clear that $\Ko=\Ass$, the associative \nso\ in the monoidal cartesian category of spaces.

The  linear embedding $\epsilon$ induces a morphism of operads
$$
\epsilon_\#\colon\calK_1\longrightarrow\Kd
$$
which sends a configuration on the line to its image under $\epsilon$ in $\BR^d$. 
This restricts to a morphism 
$$
\epsilon_\#\colon\Ko=\Ass\longrightarrow\Kd
$$
which endows $\Kd$ with the structure of a multiplicative operad.

One can thus associate a cosimplicial space $\Kd^\bullet$, which we will call \emph{Sinha's cosimplicial space}, to $\Kd$, as outlined in \refS{cos-multop}. 
In more detail, this is a cosimplicial space
$$
\xymatrix{
\Kd^\bullet = 
\big(\Kd(0) \ar@<0.7ex>[r] \ar@<-0.7ex>[r] &
\Kd(1) \ar[l] \ar[r] \ar@<1.2ex>[r]  \ar@<-1.2ex>[r] &
\Kd(2) \ar@<0.6ex>[l]  \ar@<-0.6ex>[l]
\cdots \big),
}
$$
where $\Kd(n)$ has the homotopy type of the space of configurations of $n$ points in $\BR^d$. Cofaces $d^i$ corresponds to 
``doubling'' the $i$th point of
the configuration ``infinitesimally'' in the direction given by $\epsilon$ and codegeneracies $s^i$ forget the $i$th point in the configuration.  This is explained in detail in \cite{S:TSK,S:OKS}.

We have the following important result due to Sinha. 
\begin{thm}[\cite{S:OKS}, Corollary 1.2]
For $d\geq 4$, the space of long knots modulo immersion is weakly equivalent to the homotopy totalization
of the cosimplicial space associated to the Kontsevich operad, i.e. 
$$\bEmb(\BR,\BR^d)\simeq\hoTot(\Kd^\bullet).$$
\end{thm}


\subsection{The Fulton-MacPherson operad}\label{S:Fd}


We recall here the Fulton-MacPherson operad and its relation to the Kontsevich operad. We will need this operad to establish certain formality results. Our main reference is \cite{S:comp}, although 
this operad is also studied in
\cite{AxSi,Gaiffi,GJ:OHAII,K:OMDQ,Markl:realconf}.

The Fulton-MacPherson operad, $\Fd=\{\Fd(n)\}_{n\geq0}$, is a topological operad whose $n$th term is also a compactification of $C(n,\BR^d)$
 modulo scaling and translation.  The difference between this and the compactification defining the Kontevich operad is that no identification of the collinear configurations takes place. More precisely, $\Fd(n)$ is defined as the closure of the image of the map
\begin{equation}\label{E:defFd}
(\alpha_*,\beta_*)\colon C(n,\BR^d)\to \prod_{1\leq i<j\leq n}S^{d-1}\times \prod_{1\leq i<j<k\leq n}[0,\infty]
\end{equation}
where $\alpha_*=(\alpha_{ij})_{1\leq i<j\leq n}$ is as in \refZ{E:defKd} and  $\beta_*=(\beta_{ijk})_{1\leq i<j<k\leq n}$
is defined by $\beta_{ijk}(x_1,\cdots,x_n)=\|x_i-x_j\|/\|x_i-x_k\|$ \cite[Definition 4.11]{S:comp}.  
There is a  morphism of operads
\begin{equation}\label{E:OperadEquivalence}
q\colon\Fd\longrightarrow\Kd
\end{equation}
induced by the obvious projection between the target spaces of maps \refZ{E:defFd} and \refZ{E:defKd}.  Each $q\colon\Fd(n)\to\Kd(n)$
is a homotopy equivalence \cite[Corollary 5.9]{S:comp} (see also \cite[Theorem 4.2]{S:OKS}), as is the map of operads \eqref{E:OperadEquivalence} \cite{Sal} (see also \cite[Section 2]{LV:FLDO-nov06}).
 
Denote by $\Fo(n)$ the path component in $\calF_1(n)$ containing the linearly ordered configuration $(1,\ldots,n)$ on the line.
This defines a \nso\ $\Fo$ which is homeomorphic to the Stasheff operad.  In particular  $\Fo(n)$ is the $n$th
associahedron which is a convex polytope of dimension $n-2$ (or $0$ for $n<2$) \cite[Section 4.4]{S:comp}.
The  linear embedding $\epsilon$ 
induces a morphism 
$\epsilon_\#\colon\Fo\to\Fd$
and we have a commutative diagram of \nso s
\begin{equation}\label{E:equFdKd}
\xymatrix{\Fo\ar[r]^{\epsilon_\#}\ar[d]_{q}^{\simeq}&\Fd\ar[d]_{q}^{\simeq}\\
\Ass=\Ko\ar[r]^-{\epsilon_\#}&\Kd}
\end{equation}
where the vertical maps are homotopy equivalences.

The operad $\Fd$ is not multiplicative\ifn{This used to say ``not genuinely multiplicative".  There is way too much ``genuine" everywhere which is redundant.  If we say the operad is not multiplicative, than that means what it means.  There is no ``genuinely" or ``ingenuinely" multiplicative things.  When needed, we say ``multiplicative up to homotopy", but when something is not multiplicative, it's just ``not multiplicative".} but 
since each $\Fo(n)$ is contractible, $\epsilon_\#$ in a sense endows $\Fd$ with a multiplicative structure 
``up to homotopy''.  Thus we cannot apply the Gerstenhaber-Voronov/McClure-Smith construction directly
to obtain a cosimplicial space out of $\Fd$, but we will see in \refS{Fanicnso} that one can associate a more general 
diagram to $\Fd$ which generalizes this construction. 


\section{Formality and collapse of the homology Bousfield-Kan spectral sequence}\label{S:formalityBKSS}


In this section we recall the notion of a formal diagram in rational homotopy theory and we show that the homology Bousfield-Kan  spectral
sequence of a formal cosimplicial space collapses at the $E^2$ page. We also show that formality of a finite diagram is preserved by passing
to its cosimplicial replacement.
 
We first review some classical notions in rational homotopy theory for which \cite{FHT:RHT} is the standard reference. 
Let $\BK$ be a field of characteristic $0$ and denote by $\CDGA$ the category of commutative differential graded algebras.
Let 
$$
\Apl(\,-\,;\BK)\colon\Top\to\CDGA
$$ 
be Sullivan's functor of
piecewise polynomial forms as defined in \cite[\S 10 (c)]{FHT:RHT} and recall that, for a space $X$,
$H(\Apl(X;\BK))\cong H^*(X;\BK)$. In fact $\Apl(\,-\,;\BK)$ is naturally connected by a zig-zag of quasi-isomorphisms
of cochain complexes to the singular cochains $\C^*(-;\BK)$ \cite[\S 10(e)]{FHT:RHT}.
The definition of a \emph{formal space} was introduced in \cite{DGMS}, and it can be generalized to diagrams as follows.

\begin{definition}\label{D:formal}
Let $I$ be a small category and let $\BK$ be a field of characteristic $0$. 
A functor $A\colon I\to\CDGA$ is called \emph{formal} if it is connected by a zig-zag of natural quasi-isomorphism
to its homology $\Ho(A)$. 
A diagram of spaces $F\colon
I\to\Top$ is called \emph{$\BK$-formal} if the contravariant functor $\Apl(F;\BK)\colon I\to\CDGA$ is formal.
\end{definition}

As a special case, a cosimplicial space $X^\bullet$ is formal if the diagram $X^\bullet\colon\Delta\to\Top$ is formal.
Our main interest in this is the following collapsing result.

\begin{prop}\label{P:formalcollapse}
Let $X^\bullet$ be a cosimplicial space and let
$\BK$ be a field of characteristic~$0$.
If the  cosimplicial space $X^\bullet$ is
$\BK$-formal then  the homology Bousfield-Kan  spectral sequence for $X^\bullet$ with coefficients in any field of characteristic $0$
collapses at the $E^2$ page.
\end{prop}
\begin{proof}
We consider the homology  Bousfield-Kan spectral sequence as constructed in \cite[Section 2.1]{B:HSS}.
For a cosimplicial chain complex $V_*^\bullet$,  denote by
$\{E^r(V_*^\bullet)\}_{r\geq 0}$ the spectral sequence induced by
the filtration by cosimplicial degree in the associated total complex of the bicomplex
$\prod_{m\geq0}N^m(V_*^\bullet)$.  Here $N^*$ is the conormalization as defined in
\cite[Section 2]{B:HSS} or in \cite[Chapter 8]{Weibel}.
It is clear that $\{E^r(C_*(X^\bullet;\BK))\}_{r\geq0}$
 coincides from the $E^2$ page with the homology Bousfield-Kan  spectral sequence of $X^\bullet$. 
The $\BK$-formality of the cosimplicial space and the natural equivalence between $\Apl$ and singular cochains
imply that there is a zig-zag of quasi-isomorphisms
$$\xymatrix{C_*(X^\bullet;\BK)&\ar[l]_-\simeq\cdots\ar[r]^-\simeq&H_*(X^\bullet;\BK).}$$
Therefore the spectral sequence $\{E^r(C_*(X^\bullet;\BK))\}_{r\geq
0}$ coincides from the $E^1$ page with $\{E^r(H_*(X^\bullet;\BK))\}_{r\geq 0}$. But
the latter spectral sequence collapses at $E^2$ because each $H_*(X^n;\BK)$ is a chain complex with $0$ differential 
so that the vertical differential  in the associated bicomplex is trivial. This proves the
statement for $\BK$. We also have an isomorphism of spectral sequences
$$\{E^r(C_*(X^\bullet;\BK))\}_{r\geq 0}\cong \{E^r(C_*(X^\bullet;\BQ))\otimes_\BQ\BK\}_{r\geq 0}, $$
and therefore the spectral sequence with rational coefficients
collapses at the same page. Applying this argument in the
opposite direction proves the result for any field of
characteristic $0$.
\end{proof}

If we could prove that Sinha's cosimplicial space $\calKd^\bullet$ was formal we would deduce 
immediately from the last proposition the collapse of the associated homology Bousfield-Kan  spectral sequence as claimed in \refT{BKSinha-coll}.
Unfortunately, it seems difficult to prove formality of  $\calKd^\bullet$ directly, so we will prove formality of some other diagrams approximating $\calKd^\bullet$ instead. 
But in order to deduce collapsing results from formality of these auxiliary  diagrams, we now need to show that the cosimplicial 
replacement of a formal finite diagram is a formal cosimplicial space.

\vskip 6pt
\noindent
A category $I$ is said to be \emph{finite} if it has  a finite
number of morphisms (hence of objects). 

\begin{prop}\label{P:formalcosimplrepl}
Let $I$ be a finite category, let $\BK$ be a field of
characteristic $0$, and let $F\colon I\to\Top$ be a diagram. If
$F$ is $\BK$-formal then so is its cosimplicial replacement $\Pi^\bullet F$.
\end{prop}
\begin{proof}
Recall from \cite[XI.5]{BK:LNM} and from the end of \refS{cos-hotot}  the definition of the cosimplicial replacement $\Pi^\bullet F$.
Inspired by this construction we associate to a \emph{contravariant} functor
$A\colon I\to\CDGA$ a simplicial CDGA, which we denote by $\otimes_\bullet A$
 and define as follows. Let   $\Nerve_\bullet(I)$ be  the nerve of $I$ and denote by 
$$u=\left(i_0\stackrel{f_1}{\leftarrow}\cdots
\stackrel{f_n}{{\leftarrow}}i_n \right)$$
a typical 
element of $\Nerve_n(I)$.
The  $n$th term of $\otimes_\bullet A$ 
is given by
$$
\left(\otimes_\bullet A\right)_n:=\bigotimes_{u\in
\Nerve_n(I)}A(i_0).
$$
Before defining the cofaces, notice that for $u\in \Nerve_n(I)$, taking all possible factorizations $f_i=f''\circ f'$ for a fixed $1\leq i\leq n$, 
postcomposing by all possible maps $g\colon i_{-1}\leftarrow i_0$ with source $i_0$, or  precomposing with all possible maps $h\colon i_n \leftarrow i_{n+1}$ with target $i_n$
 gives isomorphisms
$$
(\otimes_\bullet A)_{n+1}\cong\left\{\begin{array}{l}
\otimes_{u\in \Nerve_n(I)} \otimes_{f_i=f'\circ f''}A(i_{0})\\
\otimes_{u\in \Nerve_n(I)} \otimes_{g\colon i_{-1}\leftarrow i_0}A(i_{-1})\\
\otimes_{u\in \Nerve_n(I)} \otimes_{h\colon i_{n} \leftarrow i_{n+1}}A(i_0).
\end{array}
\right.
$$

The faces $\partial_i$ are defined using these isomorphisms and the multiplicative structure:
\begin{eqnarray*}
\partial_0&=&\bigotimes_{u\in
\Nerve_n(I)}\left(\bigotimes_{g\colon i_{-1}\leftarrow i_0
}A(i_{-1})\stackrel{\Psi_u}{\longrightarrow}A(i_0)\right)
\\
\partial_i&=&\bigotimes_{u\in\Nerve_n(I)}\left(\bigotimes_{f_i=f'\circ f''}A(i_0)\stackrel{\mult}{\longrightarrow}A(i_0)\right),
\,\,\,\,\,1\leq i\leq n,
\\
\partial_{n+1}&=& \bigotimes_{u\in N_n(I)}
\left(\bigotimes_{h\colon i_{n}\leftarrow i_{n+1} }A(i_{0})\stackrel{\mult}{\longrightarrow}A(i_0)\right),
\end{eqnarray*}
where $\Psi_u(\otimes_{g}x_g)$ is the product of the $A(g)(x_g)$ over the finite set of
maps $g\colon i_{-1}\leftarrow i_{0} $ with source $i_0$, for $x_g\in A(i_{-1})$.
To define the codegeneracies $\sigma_j$ we use the unit map $\eta\colon\BK\to A(i)$ associated to any CDGA and we set
\begin{eqnarray*}
\sigma_j\colon \bigotimes_{u'\in\Nerve_{n-1}(I)}A(i'_{0})\cong \bigotimes_{u\in
\Nerve_{n}(I),f_j=\id}A(i_0)\stackrel{\id\otimes\eta}{\longrightarrow} \bigotimes_{u\in\Nerve_n(I)}A(i_0).
\end{eqnarray*}
We call $\otimes_\bullet A$ the
\emph{simplicial replacement} of the contravariant functor $A$.
Notice that we need the nerve of $I$ to be finite in each degree  for the
above maps to be well defined.

When $F\colon I\to\Top$ is a formal diagram of spaces,
applying this simplicial replacement to each term in the
zigzag of CDGA quasi-isomorphisms connecting $\Apl(F)$ and
$\Ho^*(F)$ gives a zig-zag of quasi-isomorphisms of simplicial CDGAs
\begin{equation}\label{equ-quisosimplCDGA}
\xymatrix{\otimes_\bullet\Apl(F)&\ar[l]_-\simeq\cdots\ar[r]^-\simeq&\otimes_\bullet
\Ho^*(F).}
\end{equation}

Kunneth (quasi-)isomorphisms
induce a quasi-isomorphism of simplicial CDGAs
$$
\otimes_\bullet\Apl(F)\stackrel{\simeq}{\longrightarrow}\Apl(\Pi^\bullet   F),
$$
and similarly the standard Kunneth isomorphism gives an isomorphism
$$
\otimes_\bullet\Ho^*(F)\stackrel{\cong}{\longrightarrow} \Ho^*(\Pi^\bullet   F).
$$
Combining this with the quasi-isomorphisms $(\ref{equ-quisosimplCDGA})$
 proves the formality of $\Pi^\bullet  F$.
\end{proof}


\section{Convergence of Bousfield-Kan spectral sequences for cosimplicial spaces and their truncations}\label{S:SStotn}


Throughout this section, $X^\bullet$ is a simply-connected pointed cosimplicial  space. We will establish a condition on such a cosimplicial space which
guarantees that the Bousfield-Kan spectral sequences of its truncations converge. We also prove a folklore theorem which says that a left cofinal change of shape preserves the $E^2$ pages of the homotopy and homology Bousfield-Kan spectral sequences.

Our main convergence hypotheses will be the following
\begin{definition}
We say that $X^\bullet$ is \emph{well above the diagonal at page $E^r$} if  that page 
of the second-quadrant
homotopy Bousfield-Kan spectral sequence associated to $X^\bullet$ satisfies the  conditions
\begin{itemize}
\item[(i)] $E^r_{-p,q} = 0$ for $q\leq p$, and
\item[(ii)] for each $i$ there are finitely many $p$ such that $E^r_{-p,p+i}\not=0$.
\end{itemize}
\end{definition}
In other words, these conditions say that there are no terms on the main diagonal $E^r_{-p,p}$ or below, and that there are finitely many terms on each diagonal
line of slope $-1$.  Such conditions already appear in \cite{B:HSS}.

Consider now the $n$th truncation $X_{[n]}$ of $X^\bullet$. 
To prove a convergence result for the Bousfield-Kan spectral sequences of its
cosimplicial replacement $\Pi^\bullet  (X_{[n]})$, we need the following lemma relating
the homotopy spectral sequences associated to $X^\bullet$ and $\Pi^\bullet  (X_{[n]})$.

Recall the conormalisation functor $\norm^*\colon\Ab^\Delta\to\Ch^{\geq0}$ from cosimplicial 
abelian groups to non-negatively graded cochain complexes as defined in \cite[Section 8.4]{Weibel} or
in \cite[Section 2]{B:HSS}. For a cochain complex $C^*$ define its \emph{$n$th truncation} $\tau^nC^*$ to be the cochain complex given by
$$
\left\{
\begin{array}{ccl}
(\tau^nC^*)^q&=&C^q\textrm{ if }q\leq n\\
(\tau^nC^*)^q&=&0\textrm{ if }q> n
\end{array}
\right.
$$
with the differential induced in the obvious way from the differential on $C^*$.
\begin{lemma}\label{L:E2trunc}
The $E^2$ page of the homotopy Bousfield-Kan spectral sequence associated to $\Pi^\bullet  (X_{[n]})$  is given by
$$E^2_{-p,q}=\Ho^p(\tau^n \norm^*(\pi_q(X^\bullet))).$$
\end{lemma}
\begin{proof} The following proof was given to us by P.~Bousfield.
By \cite[XI.7.1]{BK:LNM} we have an isomorphism $E^2_{-p,q}\cong\lim^p_{\Delta[n]}\pi_q(X_{[n]})$,
where $\lim^p$ is the $p$th left derived functor of $\lim^0=\lim$.
Recall for example from \cite[8.4]{Weibel} the Dold-Kan Theorem which states that the conormalisation
functor 
$
\norm^*\colon \Ab^\Delta\to \Ch^{\geq0}
$
is an equivalence of abelian categories.
It is not difficult to adapt the proof of that result to get an equivalence of abelian categories
$$
\norm_n^*\colon \Ab^{\Delta[n]}\longrightarrow \Ch^{[n]}
$$
between $n$-truncated cosimplicial abelian groups and cochain complexes concentrated in degrees $0,\ldots,n$.
Here $\norm_n^*$ is defined in the obvious way by mimicking the definition of the usual conormalisation $\norm^*$.  In particular
for a cosimplicial abelian group $A^\bullet$ we have that $\norm^*_n(A_{[n]})\cong\tau^n(\norm^*(A^\bullet))$.

Notice that, as with  the usual Dold-Kan correspondence, for a truncated cosimplicial abelian group $A$ in 
$\Ab^{\Delta[n]}$ we have an isomorphism $\lim^0_{\Delta[n]} A\cong \Ho^0(\norm^*_n(A))$.
By a universal $\delta$-functor argument, $\Ho^p$ is the $p$th left derived functor of $\Ho^0$ in $ \Ch^{[n]}$.
The equivalence of abelian categories implies therefore that  $\lim^p_{\Delta[n]} A\cong \Ho^p(\norm^*_n(A))$.
Collecting these results proves our lemma.
\end{proof}
\begin{prop}\label{P:convBK}
If $X^\bullet$ is well above the diagonal at the $E^1$ page then
\begin{itemize}
\item[(i)]  The  homotopy and homology Bousfield-Kan spectral sequences associated to $X^\bullet$
converge strongly to  $\pi_*(\hoTot X^\bullet)$ and $\Ho_*(\hoTot X^\bullet)$;
\item[(ii)] The  homotopy and homology Bousfield-Kan spectral sequences associated to $\Pi^\bullet  (X_{[n]})$
converge strongly to  $\pi_*(\hoTot^n X^\bullet)$ and $\Ho_*(\hoTot^n X^\bullet)$;
\item[(iii)] ${\displaystyle \Ho_*(\hoTot X^\bullet)\cong\lim_{\leftarrow}\Ho_*(\hoTot^n X^\bullet)}$.
\end{itemize}
\end{prop}
\begin{proof}
(i) Since $X^\bullet$ is well above the diagonal at $E^1$ it is also
well above the diagonal at $E^2$. Then the statement we want is exactly the content of the results in
\cite[Proposition IX.5.7]{BK:LNM} and \cite[Theorem 3.2]{B:HSS}.

(ii) By \refL{E2trunc} we have that $\Pi^\bullet  (X_{[n]})$ is well above the diagonal at $E^2$.  Moreover $\hoTot(\Pi^\bullet  (X_{[n]}))\simeq\hoTot^n X^\bullet$.  Thus (ii) follows from the same argument as in (i).

(iii) By the hypothesis and \refL{E2trunc} we have that the connectivity of the map of spectral sequences between the
 $E^2$ pages of the homotopy spectral sequence for $X^\bullet$ and $\Pi^\bullet  (X_{[n]})$ 
 tends to infinity with $n$. By the convergence of these homotopy spectral sequences this implies that
 the connectivity of the map $\hoTot  X^\bullet\to\hoTot^n  X^\bullet$ tends to infinity
 with $n$. Therefore the same is true for the homologies.
\end{proof}

The following is part of the content of~\cite[Corollary~7.4]{S:TSK}.

\begin{prop}[Sinha]\label{P:KdWAD}
For $d\geq4$, $\Kd^\bullet$ is well above the diagonal at the $E^1$ page.
\end{prop}

The last result of this section is a comparison theorem for Bousfield-Kan spectral sequences of cosimplicial replacements of diagrams which are connected by a left cofinal functor (see end of \refS{cos-hotot}).
We could not find a proof for this folklore result in the literature so we include one here.  The proof below is due to W. Dwyer (P. Bousfield has also given us another proof).

\begin{prop}\label{P:E2cofinal}
Let $\theta\colon I\to J$ be a functor between \emph{finite} categories and let $F\colon J\to\Top$ be a $J$-diagram of (pointed) 
spaces. If $\theta$ is left cofinal then both homotopy and rational homology Bousfield-Kan spectral sequences associated to the cosimplicial
replacements $\Pi^\bullet  F$ and $\Pi^\bullet  (\theta^*F)$ agree from the $E^2$ pages.
\end{prop}
\begin{proof}
The proof mimicks that of \cite[Proposition XI.9.2]{BK:LNM} and we follow most of the notation from there. 
Recall from the proof of \refP{formalcosimplrepl} the simplicial $\CDGA$ associated to a contravariant functor of $\CDGA$s.
The Kunneth quasi-isomorphism
$$
\otimes_\bullet\Apl(F)\stackrel{\simeq}{\longrightarrow}\Apl(\Pi^\bullet F)
$$
implies that the rational homology Bousfield-Kan spectral sequence for $\Pi^\bullet F$
coincides from the $E_2$ term with the spectral sequence of the double complex $\norm_*(\otimes_\bullet\Apl(F))$,
and we have an analogous result for  $\Pi^\bullet \theta^*(F)$.

Define a bisimplicial commutative graded algebra (CGA)
$$\otimes_{\bullet,\bullet}(\Ho^*(F),\theta)$$
with the $(n,q)$ term defined by
$$
\otimes_{n,q}(\Ho^*(F),\theta)=\otimes_{(u,v,\gamma)}\Ho^*(F(j_0))
$$
where  the tensor product is taken over
$$
u=(i_0\stackrel{\alpha_1}\leftarrow\cdots\stackrel{\alpha_n}\leftarrow i_n)\in\Nerve_n(I),\quad
v=(j_0\stackrel{\beta_1}\leftarrow\cdots\stackrel{\beta_q}\leftarrow j_q)\in\Nerve_q(J),\quad
\gamma\colon\theta(i_0)\to j_q.
$$
The faces and degeneracies in both directions are obvious generalizations of those in the proof of \refP{formalcosimplrepl}.

We can take the normalization of this bisimplicial CGA with respect to the first ($i=1$) or the second ($i=2$) simplicial degree
to get a simplicial chain complex
$$
\norm^{(i)}_*\left(\otimes_{\bullet,\bullet}(\Ho^*(F),\theta)\right).
$$

In the first direction we have
$$
\left(\norm^{(1)}_*(\otimes_{\bullet,\bullet}(\Ho^*(F),\theta))\right)_q=
\otimes_{v\in\Nerve_q(J)}\,\, \norm_*\left(\otimes_{(u,\gamma)\in\Nerve_\bullet(F\downarrow j_q)}\,\Ho^*(F(j_0))\right).
$$
Denote by $\Ho^*(F(j_0))_\bullet$ the constant cosimplicial CGA which consists of the ring  $\Ho^*(F(j_0))$ in each degree 
with the identity maps as faces and degeneracies. The category of simplicial CGAs is a simplicial model category
\cite[II.3 and II.5.2.(3)]{GJ:SHT}. So we have the simplicial CGA
$$
\Nerve_\bullet(F\downarrow j_q)\otimes\Ho^*(F(j_0))_\bullet\cong\otimes_{(u,\gamma)\in\Nerve_\bullet(F\downarrow j_q)}\Ho^*(F(j_0),
$$
and by cofinality, $\Nerve_\bullet(F\downarrow j_q)$ is weakly equivalent to the simplicial set $*_\bullet$ consisting of the singleton in each degree.
Using \cite[Proposition II.3.4.]{GJ:SHT} we get a weak equivalence of simplicial CGAs 
$$
\Nerve_\bullet(F\downarrow j_q)\otimes\Ho^*(F(j_0))_\bullet\simeq*_\bullet\otimes\Ho^*(F(j_0))_\bullet=\Ho^*(F(j_0))_\bullet
$$
and deduce a weak equivalence of  simplicial chain complexes
\begin{equation}\label{E:N1bi}
\norm^{(1)}_*(\otimes_{\bullet,\bullet}(\Ho^*(F),\theta))\simeq\otimes_\bullet \Ho^*(F)
\end{equation}
where each $\Ho^*(F(j_0))$ is a chain complex concentrated in degree $0$.

In the other direction,
$$
\left(\norm^{(2)}_*(\otimes_{\bullet,\bullet}(\Ho^*(F),\theta)) \right)_n=
\otimes_{u\in\Nerve_n(I)} \norm_*\left(\otimes_\bullet\Ho^*(\theta(i_0)\downarrow F)\right),
$$
where $\theta(i_0)\downarrow F$ is the composite
$$
\theta(i_0)\downarrow J \longrightarrow J\stackrel{F}\longrightarrow\Top.
$$
Here $\theta(i_0) \downarrow J$ is the undercategory (denoted by $J\backslash\theta(i_0)$ in \cite{BK:LNM}).
This undercategory has the identity map at $\theta(i_0)$ as an initial object. It is easy to deduce, by an extra degeneracy or a spectral 
sequence argument, that $\otimes_\bullet\Ho^*(\theta(i_0)\downarrow F)$ is weakly equivalent to the constant
simplicial CGA $\Ho(F(\theta(i_0))_\bullet$. We then have a weak equivalence of  simplicial chain complexes
\begin{equation}\label{E:N2bi}
\norm^{(2)}_*(\otimes_{\bullet,\bullet}(\Ho^*(F),\theta))\simeq\otimes_\bullet \Ho^*(\theta^*F)
\end{equation}
where each $\Ho^*(\theta^*F(i_0))$ is a chain complex concentrated in degree $0$.

The left hand sides of the weak equivalences \refZ{E:N1bi} and \refZ{E:N2bi} are the $E_1$ page of  spectral sequences
computing the homology of the totalization of the double complex
obtained as the double normalization of $\otimes_{\bullet,\bullet}(\Ho^*(F),\theta)$. 
Moreover these weak equivalences show that the homology of the $E_1$ pages is concentrated on a single line.
Hence both of these spectral sequences collapse at $E_2$, and since they converge to the same thing, these $E_2$ pages are isomorphic.
We deduce that the homologies of the normalizations of $\otimes_\bullet \Ho^*(F)$ and  $\otimes_\bullet \Ho^*(\theta^*F)$
are isomorphic. Those are exactly the $E_2$ pages of the homology Bousfield Kan spectral sequences of $\Pi^\bullet F$
and  $\Pi^\bullet\theta^* F$.

The proof for the homotopy spectral sequence is similar.
\end{proof}


\section{Categories of fans, fanic diagrams, and truncated cosimplicial objects}\label{S:Fans}


In this section we introduce a sequence of finite categories $\Phi[n]$ that we call categories of  \emph{$n$-fans}.
They will serve in the next section as shapes of certain \emph{fanic diagrams} associated to morphisms of operads, generalizing the Gerstenhaber-Voronov cosimplicial object associated to a multiplicative operad. We also construct a left cofinal functor 
$\phi_n\colon\Phi[n]\to\Delta[n]$.

Recall that a \emph{planar tree} is an isotopy class of an embedding of the realization of a  contractible finite
$1$-dimensional simplicial complex in the plane.  In particular a planar tree consists 
of a finite set of \emph{vertices} and \emph{edges}. The \emph{valence} of a vertex is the number of edges ending in that vertex.
A \emph{leaf} is a vertex of valence 1. The embedding in the plane induces a clockwise cyclic order on the leaves.

\begin{definition}\label{D:fan}

\begin{itemize}
\item For a natural number $n$, an \emph{$n$-fan} is a planar tree with a  distinguished vertex called 
 the \emph{bead} such that each vertex, except maybe the bead, is of valence different from $2$ and with $n+1$ leaves other than the bead
which are labeled in the clockwise cyclic order by $0,1,\ldots,n$.  The leaf
 labeled $0$ is called the \emph{root}. The bead and vertices which are not leaves are called the \emph{non-labeled vertices}.
 
\item Define a partial order on the set of $n$-fans by declaring that $T\leq T'$ if the $n$-fan $T'$ is obtained from the $n$-fan $T$ by contracting some edges connecting
 non-labeled vertices and where the bead in $T'$ is the vertex obtained by  contracting  the subtree of $T$ containing its bead.
 
\item The category corresponding to such a poset of $n$-fans is called the \emph{$n$-fanic category} and is denoted by $\Phi[n]$.  Thus there is a unique morphism $T\to T'$ if and only if  $T\leq T'$.
 
\item A diagram shaped on the fanic category $\Phi[n]$ is called an \emph{$n$-fanic diagram}.
 
 \end{itemize}
 
\end{definition}

\begin{eg}\label{X:fans}
Figure \ref{fig0} gives examples of an $8$-fan in (a) and of two $3$-fans in (b) and (c). The bead is the vertex pictured as a small circle.
The left-hand side of Figure~\ref{fig-phi2} represents the category $\Phi[2]$.

\end{eg}

\begin{figure}[h!]
\includegraphics[width=10cm]{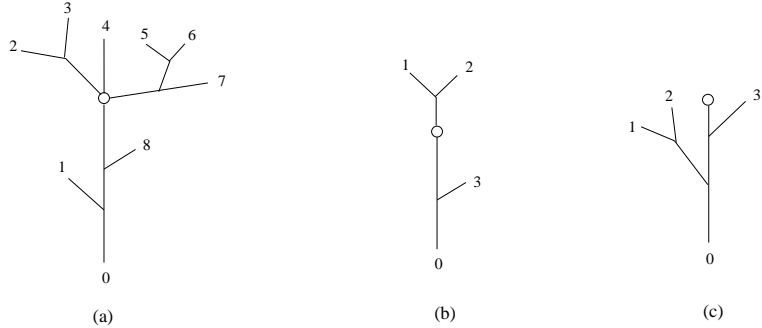}
\caption{Three examples of fans}\label{fig0}
\end{figure}

\begin{rem}
The name ``fan'' is motivated by the fact that one can define a functor from fans to rooted trees which ``opens" them along the edge emanating from the bead, much as one opens a fan.
For more details on this, see \cite{LTV:cycl-x1}. The realization of $\Phi[n]$ is homeomorphic to the barycentric subdivision of the $n$-dimensional
cyclohedron $C_n$ introduced by Bott and Taubes in~\cite{BottTaubes}.
\end{rem}

\begin{definition}\label{D:separated}
We say that a $n$-fan $T$ is \emph{$i$-separated} for  $i\in[n]$ if either
\begin{itemize}
\item the bead is not a leaf and the shortest path in $T$ joining leaves $i$ and \mbox{$(i+1)\!\mod(n+1)$} goes through the bead, or
\item the bead is a leaf and it is between leaves  $i$ and \mbox{$(i+1)\!\mod(n+1)$} in the clockwise cyclic order.
\end{itemize}
\end{definition}

\begin{eg}
In Example \ref{X:fans} the fans are $i$-separated for
\begin{itemize}
\item[(a)] $i=1,3,4,7$;
\item[(b)] $i=0,2$;
\item[(c)] $i=2$.
\end{itemize}
\end{eg}

Let  $\calPo([n])$ be the category  whose objects are non-empty subsets
of $[n]=\{0,\ldots,n\}$ and whose morphisms are inclusions. 
Define a functor 
$$\theta_n\colon\Phi[n]\longrightarrow\calPo([n])$$
by setting, for an $n$-fan $T$,
$$\theta_n(T):=\{i\in[n]:T\textrm{ is $i$-separated}\}.$$
It is immediate that $\theta_n$ is an order-preserving map between the two posets and hence a functor.

\begin{lemma}\label{L:thetacofinal}
The functor $\theta_n\colon\Phi[n]\longrightarrow\calPo([n])$ is left cofinal.
\end{lemma}
\begin{proof}
It is  easy to see  that for
any object $S\in\calPo([n])$ the overcategory $\theta_n\downarrow S$ has a terminal object and is therefore contractible.
\end{proof}

Recall the functor
${\calG}_n\colon {\calPo}([n]) \to \Delta[n]$ 
from  \cite[Definition 6.3]{S:TSK}, defined as follows.
For a non-empty subset $S\subset[n]$ consider the only order preserving bijection
$f_S\colon S\stackrel{\cong}\to[(\#S)-1]$. Define $\calG_n(S)=[(\#S)-1]$
and for a morphism $j\colon S\hookrightarrow T$ in $\calPo([n])$
let $\calG_n(j)$ be the composite $f_T\circ j\circ f_S^{-1}$. It turns out that $\calG_n$ is left cofinal
\cite[Theorem~6.7]{S:TSK}.

Now let $\phi_n$ be the composite
\begin{equation}\label{E:phin}
\phi_n:=\calG_n\theta_n\colon\Phi[n]\longrightarrow \Delta[n].
\end{equation}
\begin{thm}\label{T:phincofinal}
The functor $\phi_n:\Phi[n]\to\Delta[n]$ is left cofinal.
\end{thm}
\begin{proof}
This is immediate from \refL{thetacofinal} and the cofinality of $\calG_n$ from \cite[Theorem~6.7]{S:TSK}.
\end{proof}

\begin{eg}\label{X:phi2}
Figure \ref{fig-phi2} gives the shapes of $\Phi[2]$, $\calPo([2])$, and $\Delta[2]$ and the
functors $\theta_2$ and $\calG_2$.
\end{eg}

\begin{figure}[h!]
\begin{center}
\psfrag{Phi2}[0][0][1][0]{$\Phi[2]$}
\psfrag{Del2}[0][0][1][0]{$\Delta[2]$}
\psfrag{P}[0][0][1][0]{$\calPo([2])$}
\includegraphics[width=16cm]{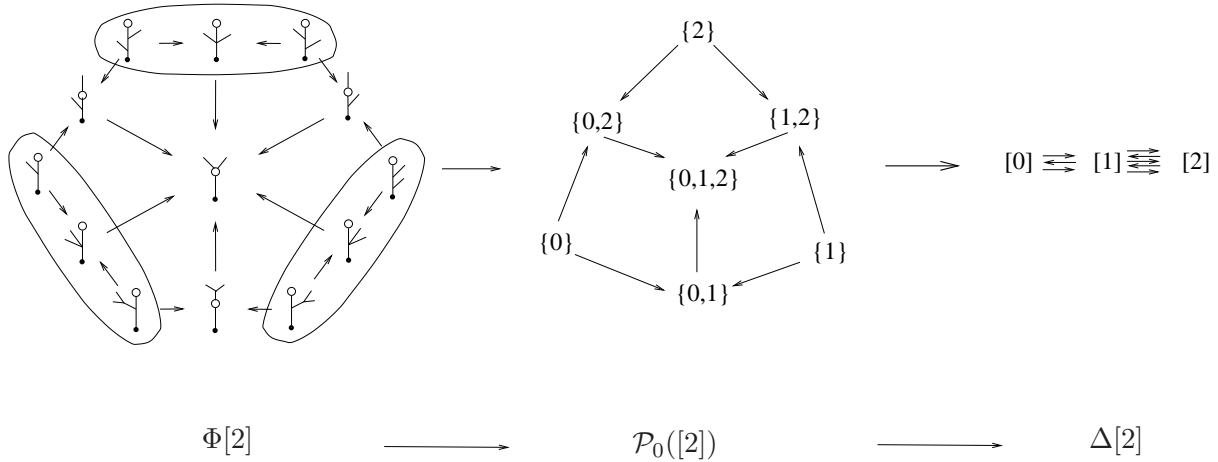}
\caption{The functor $\phi_2$}\label{fig-phi2}
\end{center}
\end{figure}


\section{The fanic diagram associated to a morphism of nonsymmetric operads}\label{S:Fanicnso}


Recall from \refS{cos-multop} the construction that associates a cosimplicial object to a multiplicative operad.  Here we generalize this construction by
associating a fanic diagram to any morphism of nonsymmetric operads.

Let 
$\mu\colon\calR\to \calM$ 
be a morphism of nonsymmetric operads in a symmetric monoidal category $\calC$.
Recall the category $\Phi[n]$ of $n$-fans defined in \refS{Fans}.
We now describe a fanic diagram
$$\fanicop\calM n\colon\Phi[n]\longrightarrow\calC$$
which will depend on $\mu$ even if only $\calM$ appears in the notation.
We first define the value of the diagram $\fanicop\calM n$ on objects. Let $T\in\Phi[n]$ be an $n$-fan. 
Orient each edge of $T$  in the unique way such that its origin is on the shortest path joining the root to the end of that edge. 
For a vertex $v$ other than the root let  $|v|$ be the number of edges emanating from $v$.
Since $v$ has exactly one incoming edge, we have $|v|=\mathrm{valence}(v)-1$.

Recall that the leaves other than the bead of the fan  $T$ are labeled and that the other vertices, including the bead,
are called the {non-labeled} vertices.
 For any non-labeled vertex $v$ of $T$ set
$$\calM(T:v):=\left\{
\begin{tabular}{ll}
$\calM(|v|)$,&if $v$ is the bead;\\
$\calR(|v|)$,&if $v$ is not the bead.
\end{tabular}
\right.
$$
Set
$$\fanicop\calM n(T):=\bigotimes_{v}\calM(T:v)$$
where the monoidal product is taken  over all the non-labeled vertices $v$ of $T$.

To define $\fanicop\calM n$ on morphisms, let $e$ be an edge of $T$ between  two non-labeled vertices $x$ and  $y$ such that $e$ is oriented from
$x$ to $y$.  Since $x$ is not the root there is a single edge that ends in $x$. Label this edge $0$ and label 
all the edges emanating from $x$ 
 with $1,\ldots,|x|$ such that they appear in the clockwise order. Hence the edge $e$ is assigned
 some label  $1\leq i\leq |x|$. Let $\overline T:=T/e$ be the fan obtained by contracting the edge  $e$ and let $\overline e$ be the vertex in $\overline T$
corresponding to that contracted edge. By definition, $\overline e$ is the \bead\ of $\overline T$ if and only if
$x$ or $y$ is the \bead\ of $T$. Notice also that the non-labeled vertices $v$ other than $x$ and $y$ in $T$ are in bijection with
 the non-labeled vertices of $\overline T$ other than $\overline e$.
The operadic structures and the morphism $\mu$  induce an obvious insertion map
$$\circ_i\colon \calM(T:x)\otimes\calM(T:y)\longrightarrow\calM(\overline T:\overline e).$$
We define the morphism
$\fanicop\calM n(T\to T/e)$ to be the composite
$$
\fanicop \calM n(T)\cong\calM(T:x)\otimes\calM(T:y)\otimes\bigotimes_{v\not=x,y}\calM(T:v)
\stackrel{\circ_i\otimes\id}{\longrightarrow}
\calM(\overline T:\overline e) \otimes\bigotimes_{v\not={\bar e}}\calM(\overline
T:v)\cong \fanicop\calM n(\overline T). 
$$

\begin{prop}\label{FanicProposition}
Let $n\geq0$ and let $\calC$ be a symmetric monoidal category.
\begin{enumerate}
\item For any morphism $\mu\colon\calR\to \calM$ of nonsymmetric operads in $\calC$ the above construction
gives a well-defined functor
$$\fanicop\calM n\colon\Phi[n]\longrightarrow\calC.$$
\item This construction is functorial, i.e. any commutative square of \nso s
\begin{equation}\label{diag-operadsquare}\xymatrix{\calR\ar[d]_f\ar[r]^\mu&\calM\ar[d]_g\\
\calR'\ar[r]^{\mu'}&\calM'}
\end{equation}
induces a natural
transformation $\widehat g\colon\fanicop\calM n\longrightarrow\fanicop{\calM'} n$.
\end{enumerate}
\end{prop}
\begin{proof}
For (1), we have already defined $\fanicop\calM n$ on objects and
contractions of one edge. Since any morphism in $\Phi[n]$ is a composition of such contractions, we only have to check that
the image under
 $\fanicop\calM n$ of a composition does not depend on the order in which we contract the edges.
 For this, it is enough to check that if $e_1$ and $e_2$ are two distinct edges between non-labeled vertices
 in an $n$-fan $T$, then
$$\fanicop\calM n(T/e_1\to T/\{e_1,e_2\})\circ\fanicop\calM n(T\to T/e_1)=
\fanicop\calM n(T/e_2\to T/\{e_1,e_2\})\circ\fanicop\calM n(T\to
T/e_2).
$$
This is an elementary check and is left to the reader.

For (2), given an object $T$ of $\Phi[n]$, set $\widehat g(T):=\displaystyle{\bigotimes_{v}}\widehat
g(T:v)$ where the tensor product is taken over the non-labeled vertices $v$ of $T$ and 
$$\widehat g(T:v):=\left\{
\begin{tabular}{ll}
$g(|v|)$,&if $v$ is the bead;\\
$f(|v|)$,&if $v$ is not the bead.
\end{tabular}
\right.
$$
This defines the desired natural transformation.
\end{proof}

Suppose moreover that the category $\calC$ is equipped with a certain class of morphisms called
\emph{weak equivalences}. This induces a class of weak equivalences on diagrams in $\calC$ by declaring
that a natural transformation between two diagrams is a weak equivalence if the map associated
to each object of the indexing category is a weak equivalence. In particular this induces a class of weak equivalences of operads in $\calC$.

\begin{prop}\label{P:wefanic}
Suppose that the symmetric monoidal category is equipped with a
class of morphisms called weak equivalences that is stable
under $\otimes$ and contains all isomorphisms. If the morphisms of
nonsymmetric operads $f$ and $g$ in the commutative square
$(\ref{diag-operadsquare})$  are weak equivalences then the
natural transformation $\widehat g$ is also a weak equivalence.
\end{prop}
\begin{proof}
Follows from the definition of $\widehat g$ and the fact that
the class of weak equivalences is stable under tensor product and contains the reordering
isomorphisms.
\end{proof}

\vskip 6pt
\noindent
Lastly, we explain in which sense the above construction is a generalization of the cosimplicial object associated 
to a multiplicative nonsymmetric operad. Recall the  functor
$\phi_n\colon\Phi[n]\to\Delta[n]$ from \refE{phin} which is left cofinal by \refT{phincofinal}.
The proof of the following is straightforward.
\begin{thm}\label{T:fanic-cosimplicial}
Let $\mu\colon\Ass\to\calM$ be a multiplicative nonsymmetric operad. Let $\calM_{[n]}\colon\Delta[n]\to\calC$ be the $n$th truncation
of the associated cosimplicial object $\calM^\bullet$ and let  $\fanicop{\calM}{n}$ be the $n$-fanic diagram asscoiated to the morphism $\mu$.
Then the following diagram commutes:
$$\xymatrix{\Phi[n]\ar[d]_{\phi_n}\ar[r]^{\fanicop{\calM}{n}}&\calC\\
\Delta[n]\ar[ru]_{\calM_{[n]}}.&}$$
\end{thm}

\begin{cor}\label{C:holimfanic}
Let $\calM$ be a multiplicative operad of spaces, let $\calM^\bullet$ be the associated cosimplicial space
and let $\fanicop\calM n$ be the associated $n$-fanic diagram for $n\geq0$. Then there is a homotopy equivalence 
$$
\holim_{\Phi[n]}\fanicop{\calM}{n}\simeq\hoTot^n(\calM^\bullet).
$$
\end{cor}
\begin{proof}
This is a consequence of \refT{fanic-cosimplicial} and \refT{phincofinal}.
\end{proof}

\begin{rem} 
Notice that all the results of this section can
be easily dualized for cooperads. In particular, to a morphism of
cooperads $\mu^*\colon\calM^*\to\calR^*$ in $\calC$ we can
associate a \emph{cofanic} diagram
$\fanicop{\calM^*}{n}\colon\Phi[n]^{op}\to \calC$
where $(-)^{op}$ means the opposite category.
We leave it as an exercise for the reader to state and prove these dual results.
\end{rem}
\begin{rem}
Notice that the map $\mu\colon\calR\to\calM$ endows $\calM$ with the structure of a bimodule over the \nso\
$\calR$. It is easy to see  that the entire discussion above can actually be applied to any such bimodule.  In particular, 
one has an associated fanic diagram for a bimodule (this motivates the notation $\fanicop\calM n$ where $\calM$ denotes that bimodule).
\end{rem} 
\begin{rem}
It is possible to describe an infinite category $\Phi$ filtered by the categories $\Phi[n]$.
One can also construct an $\infty$-fanic diagram $\fanicop\calM\infty\colon\Phi\to\calC$
whose restriction to $\Phi[n]$ is $\fanicop\calM n$, and a left cofinal functor
$\phi\colon\Phi\to\Delta$. In the case of a multiplicative operad $\calM$, we have
 $\phi^*\calM^\bullet\cong \fanicop\calM\infty$. This seems to be a more natural generalization
of the Gerstenhaber-Voronov construction since it relates directly to the cosimplicial object $\calM^\bullet$ instead of its truncations. However, we do not use this construction here since the homology Bousfield-Kan spectral sequence of the cosimplicial replacement of the infinite category $\Phi$ might be troublesome.
\end{rem}


\section{Formality of the Fulton-McPherson fanic diagram}\label{S:formalfanicFd}


Recall from \refS{Fd} the Fulton-MacPherson operad $\Fd$ which consists of suitable compactifications of configuration spaces in $\BR^d$.  A linear inclusion $\epsilon\colon\BR\to\BR^d$ induces a morphism of \nso s
$$\epsilon_\#\colon\Fo\to\Fd$$
where $\Fo$ is the principal path component of the Fulton-MacPherson operad in dimension $1$.  
From \refS{Fanicnso}, we then have an associated \emph{$n$-fanic Fulton-MacPherson diagram}
$$\fanicop\Fd n\colon\Phi[n]\longrightarrow\Top.$$

The goal of this section is to establish the following
\begin{thm}\label{T:FormalityFanic}
For $d\geq 3$ and $n\geq0$, the Fulton-McPherson \fanic\ diagram $\fanicop\Fd n$ is $\BR$-formal.
\end{thm}

The proof of this theorem is based on Kontsevich's theorem on the formality of the little $d$-disks operad, proved in  \cite[Section 3]{K:OMDQ}, with some relevant parts appearing in \cite[Appendix 8]{KS:DAODC}.  A more detailed proof can be found in \cite{LV:FLDO-nov06}.  
We now recall the main ingredients and ideas of this proof.

To prove the formality of the little $d$-disks operad, Kontsevich proves the formality of $\calF_d$, the homotopy equivalent Fulton-MacPherson operad in dimension $d$  (even if in \cite{K:OMDQ} he defines what is now called the Kontsevich operad). To do so he starts
by constructing a combinatorial 
cooperad of CDGAs as follows.  

Consider the set of finite oriented graphs with
 $n$ \emph{external} vertices (labeled from $1$ to $n$) and some other, \emph{internal}, 
vertices.  Each internal vertex is at least trivalent and is connected by a path to some external 
vertex.  No double edges or loops are allowed and an ordering of the internal vertices and edges is imposed. Such graphs are called \emph{admissible}  \cite[Definition 13]{K:OMDQ}. Denote by $\Dd(n)$ the real
vector space generated by admissible graphs with $n$ external vertices, and with certain identifications (with appropriate signs) having to do with reordering of internal vertices or edges
or reversal of orientations of edges. The degree of an admissible graph $\Gamma$ is $\deg(\Gamma):=e(d-1)-qd$
where $e$ is the number of edges and $q$ is the number
of internal vertices. A degree $+1$  differential $\operatorname{d}(\Gamma)$ on  $\Dd(n)$ is defined as the
alternating sum of the graphs obtained from $\Gamma$
by contracting an edge whose at least one vertex is internal. There is also a multiplication on $\Dd(n)$ where the product of
two admissible graphs is obtained by gluing the graphs along their common external vertices.
This equips $\Dd(n)$ with the structure of a CDGA  over $\BR$.
Moreover, the sequence $\Dd=\{\Dd(n)\}_{n\geq0}$ admits a structure of a cooperad in CDGA.

One other important ingredient, developed in \cite[appendix 8]{KS:DAODC}, is the functor of \emph{semi-algebraic differential forms}
$$\OmPA^*\colon\{\textrm{semi-algebraic sets}\}\longrightarrow\CDGA,$$
mimicking the de Rham functor of smooth differential forms on smooth manifolds, where
a semi-algebraic set is a subset of $\BR^n$ defined by  finite sets of polynomial inequalities and boolean operations.
This functor $\OmPA^*$ is naturally quasi-isomorphic to the functor $\Apl(-;\BR)$.
One fact which will be important to us is that if $X$ is a semi-algebraic set of dimension
$\dim(X)\leq m$ then $\OmPA^i(X)=0$ for $i>m$. 

Spaces $\calF_d(n)$ are semi-algebraic manifolds and Konsevich's idea is to assign to each
 admissible graph $\Gamma\in\Dd(n)$ of degree $r$ a certain 
semi-algebraic differential form $I(\Gamma)=\omega_\Gamma\in\OmPA^r(\Fd(n))$. This defines a CDGA morphism
$$I_n\colon \Dd(n)\longrightarrow\OmPA^*(\Fd(n)).$$
On the other hand, it is easy to construct an explicit CDGA morphism
$$\bar I_n\colon \Dd(n)\longrightarrow\Ho^*(\Fd(n);\BR)$$
that sends to $0$ any  admissible graphs with at least one internal vertex  and that sends
the admissible graph with a unique edge joining the external vertices $i$ and $j$, for $1\leq i<j\leq n$,
 to the generator $g_{ij}$ in the usual presentation of the cohomology of
the configuration space $$\Ho^*(\Fd(n);\BR)\cong\bigwedge(\{g_{ij}:1\leq i<j\leq n\})/\sim.$$
The exact presentation for the three equivalence relations $\sim$ (one of which is the Arnold, or three-term, relation) can be found, for example, in \cite[Section 7]{S:TSK}.

It can be proved that $\bar I_n$ is a quasi-isomorphism and, since $I_n$ is surjective on the indecomposables in cohomology, it is 
also a quasi-isomorphism. Finally, it would be nice if these were quasi-isomorphisms of cooperads, but they are not because the 
contravariant functor $\OmPA^*$ is not monoidal and so $\OmPA^*(\Fd)$ does not inherit the structure of  a cooperad. 
Nevertheless $I=\{I_n\}_{n\geq0}$
is  almost a quasi-isomorphism of cooperads.  More precisely, we have the following
\begin{thm}[\cite{LV:FLDO-nov06}]\label{T:AlmostCooperad}  There is a commutative diagram

\begin{equation}\label{diag-almostformality}
\xymatrix{
 \Dd(n) \ar[dd]_-{I_n}^-{\simeq} \ar[r]^-{\phi} & \Dd(k)\otimes \Dd({n_1})\otimes \cdots\otimes \Dd({n_k})
 \ar[d]^-{I_k\otimes I_{n_1}\otimes\cdots\otimes I_{n_k}}_-\simeq  \\
    &    \OmPA^*(\Fd(k))\otimes \OmPA^*(\Fd(n_1))\otimes\cdots\otimes\OmPA^*(\Fd(n_k))
\ar[d]_-{\simeq}^-{\textrm{Kunneth}} \\
 \OmPA^*(\Fd(n)) \ar[r]^-{\OmPA^*(\mu)} & \OmPA^*(\Fd(k)\times \Fd(n_1)\times \cdots \times \Fd(n_k))
}
\end{equation}
where $\mu$ and $\phi$ are the (co)operadic structure maps on $\Fd$ and $\Dd$ and the vertical maps are quasi-isomorphisms.
\end{thm}

On the other hand, $\{\bar I_n\colon\Dd(n)\quism\Ho^*(\Fd(n);\BR)\}_{n\geq0}$ is a quasi-isomorphism of cooperads. 
This, combined with \refT{AlmostCooperad}, is what we mean when we say that the 
Fulton-MacPherson operad is formal over $\BR$. Notice that if we work dually at the level of chains, we can use a monoidal functor to get genuine formality in the category of chain complexes.

In order to prove the formality of the Fulton-MacPherson fanic diagram, we will prove the formality of the morphism 
$\epsilon_\#\colon\Fo\to\Fd$.  The key result we need is  \refL{collinearities} below.  To explain this,
let $\COASS=\{\BR\}_{n\geq0}$ be the coassociative cooperad in $\CDGA$. 
In degree $0$, the vector space $\Dd(n)$ is generated by the single admissible graph with $n$ external vertices and no edges. 
These isomorphisms  $\Dd(n)^0\cong\BR$ induces a morphism of nonsymmetric cooperads in CDGA
$$
\epsilon\colon\Dd\longrightarrow\COASS.
$$
Let
$\eta_n\colon\BR\to\OmPA^*(\Fo(n))$ be the inclusion of constants
in degree $0$.
\begin{lemma}\label{L:collinearities}
For $d\geq 3$, diagram 
$$
\xymatrix{ 
\Ho^*(\Fd(n);\BR)\ar[d]_{\Ho^*(\epsilon_\#)}&
\ar[l]_-{\bar I}^{\simeq}\Dd(n)  \ar[r]^-{I_n}_-{\simeq}  \ar[d]^-{\epsilon_n}&
\OmPA^*(\Fd(n)) \ar[d]^-{\OmPA^*(\epsilon_\#)}
\\
\Ho^*(\Fo(n);\BR)&
\ar[l]_-{\cong} \BR  \ar[r]^-{\simeq}_-{\eta_n}&
\OmPA^*(\Fo(n))  \\
}
$$
commutes.
\end{lemma}
\begin{proof}
The left square is clearly commutative.  For the right square in case
$n\leq 1$, $\Dd(n)$ is concentrated in degree 0 and diagram clearly commutes. So assuming $n\geq 2$, let $\Gamma\in \Dd(n)$ be an admissible graph of
positive degree  with $n$ external vertices. Assume first that
$\Gamma$ has no isolated vertices, i.e. each external vertex is
an endpoint of at least one edge. Let $q$ be the number of
internal vertices of $\Gamma$ and recall that every internal
vertex has valence at least 3.  Then the number $e$ of edges of
$\Gamma$ satisfies
$$
e\geq \frac{1}{2}(n+3q).
$$
Therefore
$$
\deg(I(\Gamma))=\deg(\Gamma)=e(d-1)-qd\geq \frac{1}{2}(n+3q)(d-1)-qd=q\left(\frac{d-3}{2}\right)+n\left(\frac{d-1}{2}\right).
$$
Since $d\geq 3$ this implies that $\deg(I(\Gamma))\geq n$. On the
other hand, we have $\dim(\Fo(n))= n-2$ and so
$$
\deg(I(\Gamma))>\dim(\Fo(n)).
$$
Since for any semi-algebraic set $X$, $\OmPA^i(X)=0$ if $i>\dim(X)$, we deduce that
$$
\OmPA^*(\epsilon)(I(\Gamma))=0\in \OmPA^*(\Fo(n)).
$$

Now suppose $\Gamma$ has only $j<n$ external vertices which are not isolated.
Let $\sigma\colon \Dd(j)\to\Dd(n)$ be the composite of 
codegeneracies $\sigma_i$ inserting the $n-j$ isolated vertices. 
Hence there exists an admissible graph $\Gamma'\in\Dd(j)$ with no isolated vertex and such that
$\Gamma=\sigma(\Gamma')$. Consider 
the map $s \colon \Fd(n)\to\Fd(j)$ obtained
as the composition of $n-j$ codegeneracies $s^i$ that forget the configuration points corresponding to these isolated vertices.
 Since (co)degeneracies are induced by the operadic structure and $I$ is a map of (almost) cooperads,
 it is easy to check that the following diagram commutes:
$$
\xymatrix{
\Dd(j) \ar[r]^-{I}_-{\simeq} \ar[d]^-{\sigma} &  \OmPA^*(\Fd(j)) \ar[r]^-{\epsilon^*_\#} \ar[d]^-{s^*} &  \OmPA^*(\Fo(j)) \ar[d]^{s^*}   \\
\Dd(n)  \ar[r]^-{I}_-{\simeq}  &  \OmPA^*(\Fd({n})) \ar[r]^-{\epsilon^*_\#}   &  \OmPA^*(\Fo(n))
}
$$
Since $\Gamma'$ has no isolated vertices
it follows from the previous case that   $\epsilon_\#^*(I(\Gamma'))=0$ and so
$$
\epsilon_\#^*(I(\Gamma))=\epsilon^*_\#(I(\sigma(\Gamma')))=s^*(\epsilon^*_\#(I(\Gamma')))=0.
$$
\end{proof}

We are now ready for the proof of the formality of the Fulton-MacPherson fanic diagram.
\begin{proof}[Proof of Theorem \ref{T:FormalityFanic}]
The fanic diagram $\fanicop\Fd n$ induces a cofanic diagram
of CDGAs
$$
\OmPA^*(\fanicop\Fd n)\colon\Phi[n]^{op}\longrightarrow\CDGA.
$$
We also have a diagram of CDGA cooperads
$$\xymatrix{\Dd\ar[d]_\simeq^{\bar I}\ar[r]^{\epsilon}&\COASS\ar[d]_=\\
H^*(\Fd;\BR)\ar[r]^-{\epsilon^*}&H^*(\Fo;\BR)=\COASS}$$
which
induces quasi-isomorphic cofanic diagrams
$$
\fanicop\Dd n\simeq\fanicop{H^*(\Fd;\BR)}n\colon\Phi[n]^{op}\longrightarrow\CDGA,
$$
and of course the cofanic diagram $\fanicop{H^*(\Fd;\BR)}n$ is
$\BR$-formal.

We still need to construct a natural quasi-isomorphism
$$
\widehat I\colon \fanicop\Dd n\stackrel{\simeq}{\longrightarrow}\OmPA^*(\fanicop\Fd n).
$$
For a fan $T\in\Phi[n]$ with bead $b$, and for a non-root vertex $v$
recall that $|v|=\mathrm{valence}(v)-1$. Consider the morphism
$$
I_{|b|}\otimes\otimes_{v\not=b}\eta_{|v|}\colon\fanicop\Dd n(T)\cong\Dd(|b|)\otimes\otimes_{v\not=b}\BR\longrightarrow
\OmPA^*(\Fd(|b|))\otimes \otimes_{v\not=b}\OmPA^*(\Fo(|v|))$$
 and define
$\widehat I(T)$ as the composite of this map with the Kunneth
quasi-isomorphism
$$
\OmPA^*(\Fd(|b|))\otimes\otimes_{v\not=b}\OmPA^*(\Fo(|v|))\,\stackrel{\simeq}{\longrightarrow}\,
\OmPA\left(\Fd(|b|)\times\prod_{v\not=b}\Fo(|v|)\right).$$ We
need to show that $\widehat I$ is a natural transformation. Remember
from \refS{Fanicnso} the definition of the image of a morphism by the fanic diagram associated
to a morphism of (co)operads. Let
$e$ be an edge of $T$ emanating from $x$ and ending at $y$,
 let $i$ be the label of the output $e$
with respect to its origin $x$,
 and let
 $\overline e$ be the vertex in $\overline T=T/E$ corresponding to this contracted edge.
 Denote by $v$ any vertex of $T$ different from $x$ and $y$
or  any vertex of $\overline T$ different from $\overline e$ (these two sets of vertices are the same). We have three
cases: (i) $x=b\not=y$, (ii) $x\not= b=y$, and (iii) $x\not=b\not=y$. We will treat only case (i), the second case being
completely analoguous and the third simpler. In  case (i), $\overline e$ is the bead of $\overline T$. Consider the diagram
\begin{equation}\label{diag-formalcofanic1}
\xymatrix{
\Dd(|\overline e|)\ar[d]_-{\circ_i}\ar[rr]^-{I_{|\overline e|}}_-\simeq&
&
\OmPA^*(\Fd(|\overline e|))\ar[d]^-{\OmPA^*(\circ_i)}\\
\Dd(|x|)\otimes\Dd(|y|)\ar[d]_-{\id\otimes\epsilon}\ar[r]^-{I_{|x|}\otimes I_{|y|}}_-\simeq&
\OmPA^*(\Fd(|x|))\otimes\OmPA^*(\Fd(|y|))\ar[r]_-\simeq^-{\text{Kunneth}}\ar[d]_-{\id\otimes\epsilon^*}&
\OmPA^*(\Fd(|x|)\times\Fd(|y|))\ar[d]^-{(\id\times\epsilon)^*}\\
\Dd(|x|)\otimes\BR\ar[r]^-{I_{|x|}\otimes\eta_{|y|}}_-\simeq&
\OmPA^*(\Fd(|x|))\otimes\OmPA^*(\Fo(|y|))\ar[r]_-\simeq^-{\text{Kunneth}}&
\OmPA^*(\Fd(|x|)\times\Fo(|y|)) }
\end{equation}
The top rectangle is commutative as a consequence of the commutativity of diagram
$(\ref{diag-almostformality})$ and the definition of the insertion maps $\circ_i$ from the (co)operad structures. The
left bottom square is commutative by Lemma \ref{L:collinearities}, and the right bottom square is commutative by
naturality of the Kunneth quasi-isomorphism.

Consider now the outermost square of $(\ref{diag-formalcofanic1})$ and tensor the left side by
$\otimes_v\BR$ and the right side by $\otimes_{|v|}\OmPA^*(\Fo(|v|))$. Applying once more the Kunneth quasi-isomorphism to the right side gives
the commutative diagram
$$
\xymatrix{ \fanicop\Dd n(\overline T)\ar[r]^-{\widehat I(\overline
T)}_-\simeq\ar[d]&
\OmPA^*(\fanicop\Fd n(\overline T))\ar[d]\\
\fanicop\Dd n(T)\ar[r]^-{\widehat I( T)}_-\simeq&
\OmPA^*(\fanicop\Fd n(T)) }
$$
which is what we were after.
\end{proof}
\begin{rem}\label{R:nonmult}
Let us explain more precisely what we mean when we say that we cannot prove that the Kontsevich operad
is formal as a multiplicative operad and why we need to go through fanic diagrams. 
Lemma \ref{L:collinearities} proves that the morphism $\epsilon_\#\colon\Fo\to\Fd$ is formal,
so the weak equivalence \refZ{E:equFdKd} implies that we have a chain of quasi-isomorphisms 
of morphisms of (almost) cooperads between $\Apl(\Ass)\rightarrow\Apl(\Kd)$ and $\Ho^*(\Ass;\BR)\rightarrow \Ho^*(\Kd;\BR)$.
The two cooperads $\Apl(\Kd)$ and $H^*(\Kd;\BR)$ are comultiplicative since
$\Apl(\Ass)=\Ho^*(\Ass;\BR)=\{\BR\}_{n\geq0}$. The problem is that, in the chain of quasi-isomorphisms
joining $\Apl(\Ass)$ and $\Ho^*(\Ass;\BR)$, nothing guarantees that all the intermediate $\CDGA$s will be
coassociative even if their homologies are. 
Therefore we cannot apply the dual of the Gerstenhaber-Voronov/McClure-Smith construction to the intermediate $\CDGA$s.
One way around this would be to use ``simplicial up to homotopy'' $\CDGA$s but we wanted to avoid that. Instead, fanic diagrams allow us to work with strictly commutative diagrams.
\end{rem}


\section{Proof of the main theorem}\label{S:collHSSKd}


In this section we finally prove \refT{BKSinha-coll}.  We begin with

\begin{prop}\label{P:Kdncoll}
For $d\geq 3$ and $n\geq 0$ the rational homology Bousfield-Kan spectral sequence of the
cosimplicial replacement $\Pi^\bullet   \Kdn $ of the $n$th truncation of Sinha's cosimplicial
space $\Kd^\bullet$ collapses at the $E^2$ page.
\end{prop}
\begin{proof}
By \refT{FormalityFanic}, the Fulton-MacPherson fanic diagram $\fanicop\Fd n$ is $\BR$-formal.
By diagram \refZ{E:equFdKd} and \refP{wefanic}, we have an equivalence of diagrams
$\fanicop\Fd n\simeq\fanicop\Kd n$ and by \refT{fanic-cosimplicial} we have 
$\fanicop\Kd n=\phi_n^*(\Kdn )$ where $\phi_n\colon\Phi[n]\to\Delta[n]$
is the left cofinal functor from \refT{phincofinal}. 

Thus the finite diagram $\phi_n^*(\Kdn )$ is $\BR$-formal, and so is its cosimplicial replacement by 
\refP{formalcosimplrepl}.  By \refP{formalcollapse}, its homology spectral sequence 
collapses at the $E^2$ page. Since  $\phi_n$ is left cofinal by  \refT{phincofinal}, we deduce
from \refP{E2cofinal} that the homology spectral sequence of $\Pi^\bullet  \Kdn $ collapses at $E^2$ page.
\end{proof}

We are now ready for the proof of our main theorem. The trick will be to replace $\Kd^\bullet$ (which we do not\ifn{There was a ``don't" here, but it's bad form to use abbreviations like that in written English.} know if it is formal) 
by an associated formal cosimplicial space $\Xi^\bullet$.

\begin{proof}[Proof of \refT{BKSinha-coll}]
We first recall some ideas from rational homotopy theory.
Let $\CDGA_1$ be the full subcategory of $\CDGA$ over $\BQ$ 
of simply-connected $\CDGA$s of finite type, and let $\DGL$ 
be the category of connected differential graded Lie algebras as in \cite[\S 21 (f)]{FHT:RHT}. As 
explained in \cite[\S 22]{FHT:RHT}, there are two functors
$$
\CDGA_1\stackrel{\calL}\longrightarrow\DGL\stackrel{C^*} \longrightarrow \CDGA
$$
where $\calL(A)$ is essentially the primitive part of the cobar on the dual coalgebra $\hom(A;\BQ)$
and $C^*(L)$ is the dual of the bar construction on the enveloping algebra of $L\in\DGL$. 
The only property of these functors of interest to us is that, for $A\in\CDGA_1$, $C^*(\calL(A))$ is a Sullivan algebra quasi-isomorphic to $A$
 \cite[\S 22 (e)]{FHT:RHT}.
In other words, $C^*(\calL(-))$ can serve as a cofibrant replacement functor.

We also  have a \emph{spatial realization} functor
$$|-|\colon\CDGA\to\Top$$
defined in \cite[\S 17]{FHT:RHT}. An important  property of this functor is that a Sullivan algebra
$A$ is naturally weakly equivalent to $\Apl(|A|;\BQ)$ \cite[\S 17 (d)]{FHT:RHT}. 

Now consider the simplicial $\CDGA$ $H^*(\Kd^\bullet;\BQ)$ and define the cosimplicial space
$$\Xi^\bullet:=\left|C^*(\calL(H^*(\Kd^\bullet;\BQ)))\right|.$$
It is clear from the discussion above that this cosimplicial space is formal and has the same cohomology
as $\Kd^\bullet$.

A more surprising property, arising from the fact that the cohomology algebras of configuration spaces in
$\BR^n$ are Koszul, is that the cosimplicial space $\Xi^\bullet$ is also coformal.
This is not difficult to prove and details are in \cite[Section 4]{ALTV:coformal-x2}, in particular Corollary 4.3, 
noticing that $\Xi^\bullet$ is exactly $|\chi^\bullet|$ in that paper. 
Also explained in that paper is the consequence that there exists an isomorphism of cosimplicial groups
\begin{equation}\label{E:pichiKd}
\pi_*(\Xi^\bullet)\cong\pi_*(\Kd^\bullet)\otimes\BQ.
\end{equation}

We know by \refP{KdWAD} that $\Kd^\bullet$ is well above the diagonal at the $E^1$ page and by
 \refZ{E:pichiKd} the same is true for $\Xi^\bullet$.
Therefore the convergence results from \refP{convBK} hold for both  $\Kd$ and  $\Xi^\bullet$.

Since $\Xi^\bullet$ is formal, the same is true for $\Pi^\bullet  \Xi_{[n]}$, so the homology spectral
sequence of that cosimplicial replacement collapses at the $E^2$ page. But this $E^2$ page is clearly
isomorphic to the $E^2$ page of the homology spectral sequence of $\Pi^\bullet  \Kdn $, since
the two cosimplicial spaces have the same homology. Moreover, the second spectral sequence also 
collapses by \refP{Kdncoll}, and both spectral sequences converge by \refP{convBK}.
Therefore $\Ho_*(\hoTot^n\Xi^\bullet;\BQ)\cong \Ho_*(\hoTot^n\Kd^\bullet;\BQ)$, and by
\refP{convBK}(iii) this implies  
\begin{equation}\label{E:HTotchiKd}
\Ho_*(\hoTot\Xi^\bullet;\BQ)\cong \Ho_*(\hoTot\Kd^\bullet;\BQ).
\end{equation}

By formality of  $\Xi^\bullet$, its rational homology Bousfield-Kan spectral sequence also collapses at
 the $E^2$ page, which is isomorphic to the $E^2$ page of $\Kd^\bullet$. Since both these spectral 
 sequences converge to the isomorphic terms of \refZ{E:HTotchiKd},
 we deduce that the homology spectral sequence of $\Kd^\bullet$ also collapses at the $E^2$ page.
 \end{proof}

\bibliography{PLbib-jan07}
\bibliographystyle{plain}
\end{document}